
\input amstex

 \input amsppt.sty
 \nologo \NoBlackBoxes
 \magnification=\magstep1
 \document

\topmatter \title  The  nil Hecke ring  and singularity of Schubert
varieties \endtitle \vskip3ex
 \author SHRAWAN KUMAR\endauthor \vskip2ex
\endtopmatter
\vskip.25cm
\centerline{{\bf Introduction}}
\vskip2ex
Let $G$ be a semi-simple simply-connected  complex algebraic
group and $T \subset B$ a  maximal torus  and a Borel subgroup
respectively.  Let ${\frak h} =$ Lie $T$ be the Cartan
subalgebra of the Lie algebra Lie $G$, and  $W := N(T)/T$
 the Weyl group associated to
the pair $(G,T)$, where $N(T)$ is the normalizer of $T$ in $G$.
We can view any element $w = \overline{w}$ mod $T \in W$ as the
element (denoted by the corresponding German character) $\frak w$
of $G/B$, defined as ${\frak w} = \overline{w} B$.   For any $w \in W$, there
is associated
the Schubert variety $X_w:= \overline{Bw B/B} \subset G/B$ and
the $T-$fixed points of $X_w$  (under the canonical left action)
are precisely $I_w :=\{{\frak v} : v \in W~$ and ~ $v \leq w\}$.

We (together with B. Kostant) have defined a certain ring $Q_W(T)$
(which is the smash product of the group algebra $\Bbb Z [W]$ with
the $W-$field $Q(T)$ of rational functions on the torus $T$ ) and
certain  elements
$y_w \in Q_W (T)$ (for any $w \in W)$.  Expressing the
elements $y_w$ in the $\{\delta_v\}_{v \in W}$ basis:
$$y_w =\sum b_{w^{-1},v^{-1}} \delta_v,$$
we get the  matrix $B=(b_{w^{-1},v^{-1}})_{w,v \in W}$
with entries in $Q(T)$  (cf. Definition 2.1(d)). Analogously, we  defined the
nil Hecke ring $Q_W$ ( which is the smash product of the group
algebra $\Bbb Z [W]$ with the $W-$field $Q({\frak h})$ of rational functions on
the
Cartan subalgebra ${\frak h})$ and certain  elements $x_w \in Q_W$.
Writing
$$x_w= \sum c_{w^{-1},v^{-1}} \delta_v,$$
we get another matrix $C= (c_{w^{-1},v^{-1}})_{w,v \in W}$
with entries in $Q({\frak h})$ (cf. Definition 3.1(b)).

We prove that the formal $T$-character of the ring of functions
on the scheme theoretic tangent cone $T_{\frak v} (X_w)$~(for any
 ${\frak v}
 \in I_w)$ is nothing but  $\ast b_{w^{-1},v^{-1}}$ (cf.
Theorem 2.2), where $\ast$ is the involution of $Q(T)$ given by
$e^{\lambda} \mapsto e^{- \lambda}$.  This sharpens a result
due to Rossmann [R].  In fact this work of Rossmann, and
our own work with B. Kostant on the equivariant K-theory
of flag varieties, motivated our current work. The proof of Theorem (2.2)
requires the Demazure character formula, and occupies \S 2 of this paper. We
use this theorem to prove that $ b_{w^{-1},v^{-1}} \neq 0$ if and only if
$v\leq w$,  and in this case  it has a pole of order exactly equal to $\ell
(w)$. Similarly  $ c_{w^{-1},v^{-1}} \neq 0$ if and only if $v\leq w$ (cf.
Corollaries 3.2).

We study the graded algebra structure on the space of functions Gr ($\Cal
O_{\frak v,X_w}$) on the scheme theoretic tangent cone  $T_\frak v(X_w)$   in
\S4. Our principal result in this direction is
Theorem (4.4), which roughly asserts that the graded algebra
Gr ($\Cal O_{\frak v,X_w}$)  arises from the natural filtration of the Demazure
module $v^{-1}V_w(\lambda) $ induced from the standard filtration of the
universal enveloping algebra $U(\frak u^-) $, where
$\frak u^-$ is the nil-radical of the opposite Borel subalgebra and
$V_w(\lambda) $ is defined in \S1. We use this theorem  to derive a result due
to Carrell-Peterson asserting that for simply-laced $G$, a point $\frak v \in
X_w$ is rationally smooth if and only if the reduced tangent cone
$T^{\text{red}}_\theta(X_w)$ is an affine space for all $v \leq \theta \leq w $
(cf. Corollary 4.11).

The   principal result of our paper is a necessary and
sufficient condition for a point ${\frak v} \in X_w$ to be smooth,
 in terms of the matrix entry
$c_{w^{-1},v^{-1}}$  (cf.  Theorem 5.5 (b)). This result asserts that for any
$v\leq w \in W$,   the point $\frak{v}\in X_w$ is smooth
 $\Leftrightarrow$
 $$
 c_{w^{-1},v^{-1}}
 = (-1) ^{\ell (w)-\ell (v)}\,
 \prod_{\beta \in S(w^{-1},v^{-1})}\ \beta ^{-1} ,
 $$
 where  $S(w^{-1},v^{-1}) :=   \{\alpha \in \Delta_+ : v^{-1}r_\alpha \leq
w^{-1} \}$.

There  is a very similar criterion for a point $\frak{v}\in X_w$ to be
rationally smooth (cf. Theorem 5.5(a)). This criterion of rational smoothness
can be easily deduced  by combining some results of Dyer and Carrell-Peterson,
but we give a different geometric proof as that proof is used  crucially  to
prove our criterion of smoothness mentioned above ( i.e. Theorem 5.5(b)).

 It should be mentioned that the elements
$c_{w^{-1},v^{-1}}$ (as well as $b_{w^{-1},v^{-1}}$) are defined
combinatorially and admit closed expressions (cf. Lemma 3.4).

The nil Hecke ring approach to singularity, developed in this paper, is
applied to some specific examples discussed in \S\S6 and 7.  In \S6, we
determine the precise singular locus of any Schubert  variety in any rank-2
group (cf. Proposition 6.1). I believe this result should be well
known, but I did not find it explicitly written down in the literature. In \S7,
we use our Theorem (5.5) to study the smoothness (and rational smoothness) of
codimension one Schubert varieties $X_i $ in any $G/B$. Proposition (7.4)
(resp. Corollary 7.6) gives a criterion for a point $\frak v \in X_i$ to be
smooth (resp. rationally smooth). This criterion is applied to give a complete
list of codimension one smooth
(as well as rationally smooth) Schubert varieties in any $G/B$ (cf. Proposition
7.8).

Finally in \S8, we extend  our main result giving the criterion of smoothness
to arbitrary (not even symmetrizable)  Kac-Moody groups (cf. Theorem 8.9).  We
also extend our result determining the formal character of the ring of
functions on the scheme theoretic tangent cone at any $\frak v \in X_w$ to
arbitrary Kac-Moody groups (cf. Theorem 8.6).
The proofs in the Kac-Moody case are similar to the finite case, and hence we
have been brief and outlined only the necessary changes.

There are other criteria for smoothness due to Lakshmibai-Seshadri
(for classical groups) \cite{LS} \cite{L}, Ryan (for SL(n)) \cite{Ry}, ...; and
for rational
smoothness due to Kazhdan-Lusztig \cite{KL},  Carrell-Peterson \cite{C},
Jantzen \cite{J},  ... ; and
by works of Deodhar and Peterson  rational smoothness implies
smoothness for simply-laced groups. It may be mentioned that
 our criterion for smoothness (as in Theorem 5.5(b) )
 is applicable to all $G$ uniformly, in contrast to the above
mentioned criteria for smoothness. We refer the reader to two
survey articles, one by Carrell \cite{C2}, and the other by Deodhar
\cite{D2}.

The main results of this paper were announced in \cite{Ku2}.

 \vskip2ex
\flushpar
{\bf Acknowledgements.}  I am grateful to W. Rossmann, D.
Peterson, and J. B. Carrell for explaining to me their (then partly
unpublished) works. I also thank J. Wahl and M. Schlessinger for some helpful
conversations. This work was partially supported by the
NSF grant No. DMS-9203660.
\vskip2ex

\centerline{{\bf 1. Notation}}

\vskip1ex

For a complex vector space $V$ (possibly infinite dimensional) ,
$V^\ast$ denotes its full vector space dual. For a finite set $S , \#S$
denotes its cardinality.

Let $G$ be a semi-simple simply-connected complex algebraic
group, and let $B$ be a fixed Borel subgroup and $T \subset B$ a
maximal torus.  Let $B^-$ be the (opposite) Borel subgroup such
that $B^-\cap B=T$.  We denote by $U$ (resp. $U^-)$ the
unipotent radical of $B$ (resp. $B^-)$.  Let ${\frak g},{\frak
b},{\frak b}^-, {\frak u},{\frak u}^-, {\frak h}$ be the Lie
algebras of the groups $G,B,B^-,U,U^-,T$ respectively.  Let
$\triangle \subset {\frak h}^*$ (resp. $\triangle_+)$ denote the
set of roots for the pair $(G,T)$ (resp. $(B,T))$.  Let
$\{\alpha_1, \cdots, \alpha_n\}$ be the set of simple roots in
$\triangle_+$ and let $\{\alpha_{1}^{\vee},
\cdots,\alpha_{n}^{\vee}\}$ be the corresponding (simple) coroots
 (where $n=~ $rank $G)$.

Let $W :=N(T)/T$ be the Weyl group (where $N(T)$ is the normalizer
of $T$ in $G)$ of $G$.  Then $W$ is a Coxeter group, generated
by the simple reflections $\{r_1, \cdots, r_n\}$ (where $r_i$ is
the reflection corresponding to the simple root $\alpha_i)$.  In
particular, we can talk of the length $\ell (w)$ of any element
$w \in W$. We denote the identity element of $W$ by e.

Let ${\frak h}_{\Bbb Z }^{\ast}:=\{ \lambda \in {\frak h}^* :
\lambda (\alpha_{i}^{\vee}) \in \Bbb Z $, for all $i\}$ be the
set of integral weights and ~$D:=\{\lambda \in {\frak
h}_{\Bbb Z }^{*}: \lambda (\alpha_{i}^{\vee}) \geq 0,$ for all $i\}$
(resp. $D^o :=\{\lambda \in {\frak
h}_{\Bbb Z }^{*}: \lambda (\alpha_{i}^{\vee}) > 0,$ for all $i\}$) the set of
dominant (resp. dominant regular) integral weights.  For any $\lambda \in
D$ and $w \in W$, we denote by $V(\lambda)$ the irreducible representation of
$G$ with highest weight $\lambda$, and  $V_w (\lambda)$
is the smallest $B-$submodule of $V(\lambda)$ containing the
extremal weight vector $e_{w\lambda}$ (of weight $w
\lambda)$.  Let $R(T):=\Bbb Z  [X(T)]$ be
the group algebra of the character group $X(T)$ of the torus
$T$.  Then $\{e^{\lambda}\}_{\lambda \in {\frak
h}^{*}_{\Bbb Z }}$ are precisely the elements of $X(T)$.  Let $Q(T)$
be the quotient field of $R(T)$.  Clearly $W$ acts on $Q(T)$ and
moreover $Q(T)$ admits an involution $*$  (i.e. a field automorphism
of order 2) taking $e^{\lambda} \mapsto e^{- \lambda}$.

For any $w \in W$, the Schubert variety $X_w$ is by
definition the closure $\overline{BwB/B}$ of $BwB/B$ in $G/B$ under the Zariski
topology (where the notation $BwB/B$ means $B\overline{w}B/B$ for any
representative $\overline{w}$ of $w$ in $N(T)$ ).  Then $X_w$ is an irreducible
(projective) subvariety of
$G/B$ of dimension $\ell(w)$.
We can view any element $w=\bar{w}\mod T \in W$ as the element
(denoted by the corresponding German character) $\frak w$ of $
G/B$, defined as $\frak w=\overline{w}B$.
 By the Bruhat decomposition, any
$\frak v$ such that $v \leq w $ belongs to $ X_w$, where $\leq $ is the Bruhat
(or
Chevalley) partial order in $W$.  The Schubert variety $X_w$ is
clearly $B-$stable (in particular $T-$stable), under the left
multiplication of $B$ on $G/B$.  The $T-$fixed points of
$X_w$ are precisely $I_w:=\{{\frak v} : v \in W$ and $ v \leq w\}$.
For any variety $X$
over $\Bbb C$, we denote by $\Bbb C [X]$ the ring of global regular
functions $X\rightarrow \Bbb C$.  For any $\lambda \in
\frak h^\ast_{\Bbb Z}$,  let
$\Bbb C_\lambda$ be the 1-dimensional representation of $B$ given by the
character $e^\lambda$ and let $\frak L(\lambda)$ be the line bundle
on $
G/B$ associated to the principal $B$-bundle $G\rightarrow
G/B$ via the representation ${\Bbb C}_{-\lambda}$ of $B$.

\vskip 7mm

\centerline{{\bf 2. Character of the ring of functions on the tangent cone
of $\bold X_{\bold w}$}}

\vskip1ex

We follow the notation as in $\S$1.
\vskip1ex
\flushpar
{\bf (2.1) Definitions.} (a)~~ For any local ring $R$ with
maximal ideal ${\frak m}$, define the graded $R/{\frak m}-$algebra:
$$\text{gr}~ R:={\displaystyle{\sum_{n \geq 0}}} {\frak m}^n/ {\frak
m}^{n+1}.$$

Let $X$ be a scheme of finite type over an algebraically closed
field  and let $x$ be a closed point of $X$.  Then the {\it
tangent cone} $T_x(X)$ of $X$ at $x$ is, by definition (cf. [M,
Chapter 3, \S3]), Spec (gr ${\Cal O}_x)$, where
${\Cal O}_x = {\Cal O}_{x,X}$ is the local ring at
$x \in X$.
\vskip1ex
\flushpar
(b)~~ Let $\widetilde{R(T)}$ be the set of all the formal
sums ${\displaystyle{\sum_{e^{\lambda} \in X (T)}}}~~
n_{\lambda} e^{\lambda}$, with arbitrary $n_{\lambda} \in
\Bbb Z $ (we allow infinitely many of the $n'_{\lambda}$s to be
non-zero).  Even though $\widetilde{R(T)}$ is not a
ring, it has a canonical $R(T)-$module structure (got by the
multiplication).  We define the $Q(T)-$module
$\widetilde{Q(T)}$ as $Q(T) \otimes_{R(T)}
\widetilde{R(T)}$. Since $Q(T)$ is a flat  $R(T)-$ module,
$Q(T)$ canonically embeds in $\widetilde{Q(T)}$.
\vskip1ex
\flushpar
(c)~~ A $T-$module $M$ is said to be a {\it weight module} if
$M={\displaystyle{\oplus_{e^{\lambda} \in X (T)}}}
M_{\lambda},$ where $M_{\lambda}:= \{m \in
M:tm=e^{\lambda} (t)m\}$ is the $\lambda-$th weight space.  A
weight module $M$ is said to be an {\it admissible $T-$module}
if dim $M_{\lambda}< \infty$, for all $e^{\lambda} \in
X(T)$.

For any admissible $T-$module $M$, one can define its {\it
formal character}~ ch~$M:={\displaystyle{\sum_{e^{\lambda}
\in X(T)}}}$ (dim $M_{\lambda})~ e^{\lambda}$ as an element of
$\widetilde{R(T)}$.
\vskip1ex
\flushpar
(d) ~~ {\it The ring }$Q(T)_W  ~$([KK2, Section 2]):  Let $Q(T)_W$ be the
smash product of the $W-$field $Q(T)$ with the group algebra
$\Bbb Z [W],$ i.e., $Q(T)_W$ is a free right $Q(T)-$module with
basis $\{\delta_{w}\}_{w \in W}$ and the multiplication is
given by:

$$ (\delta_{w_1} q_1).(\delta_{w_2} q_2)= \delta_{w_1w_2}
(w_{2}^{-1} q_1)q_2,\ \text{for}\  q_1,q_2 \in Q(T) \
\text{and} \ w_1,w_2 \in W. \leqno{(1)}$$

For any simple reflection $r_i, 1 \leq i \leq n$, define the
element $y_{r_i} \in Q(T)_W$ by:
$$y_{r_i}=(\delta_e+ \delta_{r_i}) \frac{1}{(1-e^{-
\alpha_i})}.\leqno{(2)}$$
Now, for any $w \in W$, define $y_{w} \in Q (T)_W$ by
$$ y_w =y_{r_{i_1}} \cdots y_{r_{i_p}}, \leqno{(3)}$$
where $w=r_{i_1} \cdots r_{i_p}$ is a reduced decomposition.
By [KK2, Proposition 2.4], $y_w$ is well defined. Write
$$y_{w}= {\displaystyle{\sum_{v}}} b_{w^{-1},v^{-1}} \delta_v,
\leqno{(4)}$$
for some (unique) $b_{w^{-1},v^{-1}} \in Q(T)$.  It can be esily seen that
$b_{w^{-1},v^{-1}} = 0$ unless $v \leq w$ (cf. [KK2, Proposition 2.6]).

The ring $Q (T)_W$ has a canonical representation in  $Q (T)$ defined
by
$$ (\delta_w q_1).q_2= w(q_1q_2). \leqno{(5)}$$
It is easy to see that for any $r_i$ , $y_{r_i} . R(T) \subset R(T)$ , in
particular, $y_w. R(T) \subset R(T)$ for any $w\in W$.

\vskip1ex
Since ${\frak v} \in X_{w}$ is fixed under the action of $T$
(cf. \S1), the local ring ${\Cal O}_{{\frak v},X_w}$ at ${\frak v} \in
X_w$ is canonically a $T-$module.

\proclaim{ (2.2) Theorem} {\it Take any $v \leq w \in W$.  Then} {\text gr} $
{\Cal O}_{{\frak v},X_w}$ {\it is an admissible $T-$module and  moreover
$$\text{ch~(~gr}~ {\Cal O}_{{\frak v},X_w})= * b_{w^{-1},v^{-1}},$$
as elements of $\widetilde{Q(T)}$, where {\text ch} (which is
an element of $\widetilde{R(T)})$ is to be thought of as
the element} $1 \otimes\ \text{{\text ch}}$ of $\widetilde{Q(T)}:= Q(T)
\otimes_{R(T)} \widetilde{R(T)}$.

{\it In particular,} {\text ch(gr} ${\Cal O}_{{\frak v},X_w}) \in Q(T)$.
\endproclaim

Before we come the proof of Theorem (2.2), we need the following
preparation.

We recall the following very simple lemma without proof.

 \proclaim{(2.3) Lemma} Let $Y$ be an irreducible
 projective variety with an ample line bundle
 $\frak{L}$ on Y, together with a non-zero
 $\sigma\in H^0(Y,\frak{L})$.  Define the variety
 $Y^o:= Y\backslash Z(\sigma)$, where $Z(\sigma)$ is the
 zero-set of $\sigma$.  Then $Y^o$ is affine and  moreover for any
 $f\in \Bbb C[Y^o]$, there exists a $n>0$ (depending upon $f$) such that the
 section $f\cdot \sigma^n$ (of
 $H^0(Y^o,\frak{L}^{\otimes n})$) extends as an
 element of $H^0(Y,\frak{L}^{\otimes n})$.
 \endproclaim

\proclaim{ (2.4) Lemma}  {\it Given any $f \in {\Bbb C}[U^-]$, there
exists a large enough  $\lambda \in D$ (i.e. $\lambda
(\alpha_{i}^{\vee}) >>0$, for all the simple coroots
$\alpha_{i}^{\vee})$ and $\theta \in V(\lambda)^*$ such
that
$$f(g) =\langle \theta,g e_{\lambda}\rangle,\ \text{for} \ g \in U^-,$$
where $e_{\lambda}$ is a non-zero highest weight vector of $V
(\lambda)$.

Moreover, for any $v \leq w \in W,~ f$ vanishes on $(v^{-1}
\overline{BwB}) \cap U^{-} \Leftrightarrow  \theta
\in (V (\lambda) / v^{-1} V_w (\lambda))^*$.  }
\endproclaim

 \demo{Proof} The first part is due to Andersen
  and also  Cline-Parshall-Scott \cite{CPS, \S5}.  However, for
 completeness, we give a proof.

   By the Borel-Weil theorem (for any $\lambda\in D$),
 $\chi :V(\lambda)^\ast \overset\sim\to{\rightarrow}
 H^0(G/B,\ \frak{L}(\lambda))$, where for any $\phi
 \in V(\lambda)^\ast , \chi(\phi)$ is given by the section
 $\chi (\phi)(gB)= (g,g^{-1}\phi_{\vert_{\Bbb C e_{\lambda}}})$ mod $B$.
 (Observe that $\Bbb Ce_\lambda\subset V(\lambda)$ is a
 one-dimensional representation of $B$ corresponding to
 the character $e^{\lambda}$ and hence $( \Bbb Ce_\lambda
 )^\ast$ corresponds to the character $e^{-\lambda}$.)  Let $\phi_\lambda \in
V(\lambda)^\ast$ be the element defined by $\phi_\lambda (e_\lambda) =1$ and
$\phi_\lambda (v) =0$ , for any weight vector $v\in V(\lambda)$ of weight $\mu
\neq \lambda$. Consider $U^-\approx U^-\cdot \frak{e}\subset
  G/B$ as an open subset and take any (ample) line bundle $\frak L(\lambda_o)$
on $G/B$ for $\lambda_o \in D^o.$  Taking the section $\sigma = \chi
(\phi_{\lambda_o})$ of $\frak L(\lambda_o)$ and
applying Lemma (2.3), we get
 the first part of the lemma for $\lambda = n\lambda_o$ (for some $n>0$).
(Observe that $Z(\sigma)= G/B\setminus U^-.\frak e$ , since $\lambda_o$ is
regular.)

 Let $\tilde{f}:G\to \Bbb C$ be the  extension of $f$
 given  by $\tilde{f}(g)= \langle
 \theta,ge_\lambda\rangle$.  Then, since
 $v^{-1}\overline{BwB}$ is an irreducible subvariety of $G$, and
 by Bruhat decomposition $v^{-1} \overline{BwB} \cap
 U^-B$ is non-empty open subset of $v^{-1}\overline{BwB}$,
 $$\spreadlines{2\jot}\align
 f \text{ vanishes on }v^{-1}\overline{BwB} \cap U^-
 &\Leftrightarrow \tilde{f} \text{ vanishes on
 }v^{-1} \overline{BwB} \cap (U^-\cdot B) \\
 &\Leftrightarrow \tilde{f} \text{ vanishes on }
 v^{-1}\overline{BwB}\\
 &\Leftrightarrow \tilde{f} \ \text{vanishes on}\
 v^{-1}BwB\\
 &\Leftrightarrow \langle \theta,\
 v^{-1}Bwe_\lambda\rangle =0 \\
 &\Leftrightarrow \langle \theta,\ v^{-1}V_w (\lambda)
 \rangle =0  .
 \endalign
 $$
This proves the lemma. \qed
\enddemo

 For any $\lambda\in D$, define the map
 $$
 \varphi_\lambda: V(\lambda)^\ast \otimes \Bbb C_\lambda
 \rightarrow \Bbb C[U^-]
 $$
 by $ \varphi _\lambda(\theta\otimes e_\lambda) (g)=
 \langle \theta, ge_\lambda \rangle$, for $\theta\in
 V(\lambda)^\ast ,\ g\in U^-$ and $e_\lambda\in \Bbb C_\lambda$;
 where $\Bbb C_\lambda\subset V(\lambda)$ is identified
 as the highest weight space.

 \proclaim{(2.5) Lemma} $\varphi _\lambda$ is $T$-equivariant
 with respect to the adjoint action of $T$ on $U^-$, and
 is an injective map.
 \endproclaim

 \demo{Proof} For any $t\in T$,
 $$\align
 \varphi_\lambda (t\theta\otimes te_\lambda) (g) &=
 \langle t\theta,gte_\lambda \rangle \\
 &= \langle \theta,t^{-1} g t e_\lambda \rangle \\
 &= (t\cdot \varphi _\lambda (\theta\otimes e_\lambda))
 (g).
 \endalign
 $$
 This proves the $T$-equivariance of $\varphi _\lambda$.

 To prove  the injectivity of $\varphi _\lambda$, take $\theta\otimes
 e_\lambda\in \ker \varphi_\lambda $, i.e., $\langle
 \theta,ge_\lambda \rangle=0$, for all $g\in U^-$.  Hence
 $\langle \theta,gb e_\lambda \rangle =0$ for all $g\in
 U^-$ and $b\in B$.  In particular, by the density of $U^-B$
 in $G$ and the irreducibility of $V(\lambda)$, we get $\langle
 \theta ,\ V(\lambda)\rangle=0$, i.e., $\theta=0$, proving
 the injectivity of $\varphi _\lambda$.\qed
 \enddemo

\flushpar
 {\bf \S(2.6)} For any $\gamma\in D$, let us choose a highest weight vector
 $e_\gamma \in V(\gamma)$, and define (for any $\lambda, \mu
 \in D$)
 $$
 V(\lambda+\mu ) \overset {i_{\lambda,\mu
 }}\to{\hookrightarrow} V(\lambda)
 \otimes  V(\mu) \overset {\text{Id}\otimes \pi_\mu
 }\to{\longrightarrow}
 V(\lambda) \otimes  \Bbb C_{\mu },
 $$
 where $i_{\lambda,\mu }$ is the unique $G$-module map
 taking $e_{\lambda+\mu }\mapsto e_\lambda\otimes e_\mu
 $ and $\pi_\mu : V(\mu) \rightarrow\Bbb C_\mu$  is the  T-equivariant
projection onto the
 highest weight space $\Bbb C_\mu =\Bbb
 Ce_\mu \subset V(\mu)$.  We
 denote the composite map $( \text{Id}\otimes \pi_{\mu }
 )\circ i_{\lambda,\mu }: V(\lambda+\mu ) \rightarrow
 V(\lambda)\otimes \Bbb C_\mu $ by
 $\hat{\delta}_{\lambda,\mu }$.  Dualizing the above, we
 get the map
 $$\gather
 \bar{\delta}_{\lambda,\mu } : V(\lambda)^\ast  \otimes
 \Bbb C_{-\mu} \rightarrow V(\lambda+\mu )^\ast ~,\\
 \intertext{and hence the map}
 \delta_{\lambda,\mu } = \bar{\delta}_{\lambda,\mu }
 \otimes \text{Id}\, : V(\lambda)^\ast  \otimes \Bbb
 C_\lambda \approx
 V(\lambda)^\ast \otimes \Bbb C_{-\mu } \otimes \Bbb
 C_{\lambda+\mu } \rightarrow V(\lambda+\mu )^\ast
  \otimes \Bbb
 C_{\lambda+\mu }.
 \endgather
 $$
 It is easy to see that $\delta_{\lambda,\mu }$ is
 injective.  Moreover, the following diagram is
 commutative:
  $$
  \gather
  V (\lambda)^\ast  \otimes  \Bbb C_{\lambda} \overset
 {\delta_{\lambda,\mu }}\to{\hookrightarrow} \quad
  V(\lambda+\mu )^\ast \otimes \Bbb C_{\lambda+\mu} \\
   \sideset{\varphi _\lambda}\and\to\searrow
 \qquad\qquad\qquad \sideset\and\varphi _{\lambda+\mu
 }\to\swarrow\\
  \Bbb C [U^-]
  \endgather
 $$
 By virtue  of Lemma (2.4), for any $\lambda\in D$ and $v\leq w \in W$, we get
the
 injective map
 $$
 \varphi _\lambda(v,w) : (v^{-1} V_w(\lambda))^\ast  \otimes
 \Bbb C_\lambda \hookrightarrow \Bbb C[(v^{-1}\overline{BwB}) \cap U^-] ,
 $$
 by restricting the map $\varphi _\lambda$.

 \proclaim{ (2.7) Lemma}  For $\lambda,\mu \in D$ and $v\leq
 w\in W$, $\hat{\delta}_{\lambda,\mu }
 (v^{-1}V_w(\lambda+\mu )) = v^{-1}V_w(\lambda)\otimes
 \Bbb C_\mu $.  In particular, there exists a unique map
 $\delta_{\lambda,\mu }(v,w)$ making the following
 diagram commutative:
 $$\CD
 V(\lambda)^\ast  \otimes \Bbb C_\lambda  @>>>
 (v^{-1}V_w(\lambda))^\ast \otimes \Bbb C_\lambda \\
 @VV\delta_{\lambda,\mu }V   @VV\delta_{\lambda,\mu}(v,w)V
 \\
 V(\lambda+\mu )^\ast \otimes \Bbb C_{\lambda +\mu }
 @>>>  ( v^{-1} V_w(\lambda+\mu ) )^\ast \otimes \Bbb
 C_{\lambda+\mu }
 \endCD$$
 where the horizontal maps are the canonical restriction
 maps.  Moreover, $\delta_{\lambda,\mu }(v,w)$ is
 injective.
 \endproclaim

 \demo{Proof}  For $b\in B$,
$$\hat{\delta}_{\lambda,\mu }
 (\bar{v} ^{-1} b\bar{w}e_{\lambda+\mu }) =
 \bar{v}^{-1}b\bar{w} e_\lambda \otimes [
 \bar{v}^{-1}b\bar{w}e_\mu]_{\mu}, \leqno{(1)}$$
where $[x]_{\mu }$
 denotes the component of $x\in V(\mu)$ in the
 $\mu^{\text{th}}$ weight space, and $\bar{v}$ is a
 representative of $v$ in $N(T)$.  Define the closed
 subvariety $Y\subset B$ by $Y= \{ b\in B : [
 \bar{v}^{-1}b\bar{w} e_\mu  ]_{\mu }=0\}$.  Then
 $Y\neq B$, for otherwise $e_{\mu }\notin
 v^{-1}V_w(\mu)$, which is a contradiction (since $v\leq
 w$ by assumption).  Hence for $b\in B\backslash  Y$,
 $\hat{\delta}_{\lambda,\mu }(\bar{v}^{-1} b\bar{w}
 e_{\lambda+\mu}) = \bar{v}^{-1} b\bar{w}
 e_\lambda\otimes e_{\mu }$, up to a non-zero scalar.
 But since $\hat\delta_{\lambda,\mu }(v^{-1}V_w(\lambda+\mu
 ))$ is a (closed) linear subspace and $B\backslash Y$ is
 dense in $B$,
 $$
 v^{-1} V_w (\lambda) \otimes \Bbb C_{\mu } \subset
 \hat{\delta}_{\lambda,\mu }(v^{-1}V_w(\lambda+\mu )) .
 $$
 The inverse inclusion is clear from (1).  This proves
 the first part of the lemma.  The \lq in particular'
 statement follows immediately from dualizing the map
 $$
 \hat{\delta}_{{\lambda,\mu
 }_{\vert_{v^{-1}V_w(\lambda+\mu )}}} : v^{-1}
 V_w(\lambda+\mu ) \longrightarrow
 v^{-1} V_w(\lambda) \otimes \Bbb C_\mu .
 $$
 The injectivity of $\delta_{\lambda,\mu }(v,w)$ follows
 from the surjectivity of $\hat{\delta}_{{\lambda,\mu }
 _{\vert_{v^{-1}V_w(\lambda+\mu )}}}.$ \qed
\enddemo

 By virtue of the above lemma, we get the following
 commutative diagram:
 $$
  \gather
  (v^{-1} V_w(\lambda))^\ast \otimes \Bbb C_\lambda
 \quad
  \overset{\delta_{\lambda,\mu }(v,w)}\to\hookrightarrow
 \quad
  (v^{-1} V_w(\lambda+\mu ))^\ast \otimes \Bbb
 C_{\lambda+\mu} \\ \vspace{4\jot}
  \sideset{\varphi _\lambda(v,w)}\and\to\searrow
 \qquad\qquad\qquad
 \sideset\and{\varphi _{\lambda+\mu }(v,w)}\to\swarrow\\ \vspace{2\jot}
  \Bbb C [v^{-1}\overline{BwB} \cap U^-].
  \endgather
 $$
 \flushpar
 {\bf (2.8) Definition.}  Define a partial order
 $\prec$ in $D$ as follows:
 $$
 \lambda\prec \mu \Leftrightarrow \mu - \lambda \in D.
 $$
 Taking the limit of the maps $\varphi _\lambda(v,w)$,
 we get the $T$-equivariant map
 $$
 \varphi (v,w) :  \underset{\lambda\in D}\to{\text{limit}_\rightarrow}
 ((v^{-1}V_w(\lambda))^\ast \otimes \Bbb C_\lambda)
 \longrightarrow
 \Bbb C[v^{-1}\overline{BwB} \cap U^-].
 $$

 \proclaim{(2.9) Proposition} The above map $\varphi (v,w)$ is
 an isomorphism, for all $v\leq w\in W$.
 \endproclaim

 \demo{Proof} Injectivity of the map $\varphi (v,w)
 $ is clear from the injectivity of the maps
 $\varphi_\lambda(v,w)$ and $\delta_{\lambda,\mu
 }(v,w)$.  Surjectivity of $\varphi (v,w)$ follows from
 Lemma (2.4). \qed
 \enddemo

 \flushpar {\bf (2.10) Definition.} For any directed set
 $\Lambda$ and any sequence
 $
 \theta: \Lambda \rightarrow \widetilde{R(T)},
 $
 given as $\theta(\alpha)= \sum_{e^\lambda \in X(T)}
n_\lambda(\alpha)e^\lambda$ with
 $n_\lambda(\alpha)\in \Bbb Z$, we say that $\underset
 {\alpha\in \Lambda}\to{\text{limit}} ~\theta(\alpha) = \sum
 n_\lambda e^\lambda$ , if for any $e^\lambda\in X(T)$,
 there exists $\alpha_\lambda\in \Lambda$ such that
 $n_\lambda(\alpha)= n_\lambda$ for all $\alpha\geq
 \alpha_\lambda$.  Of course $\underset{\alpha\in
 \Lambda}\to{\text{limit}}  ~\theta(\alpha)$ may not exist
 in general.

Observe that if ~ $\underset{\alpha\in
 \Lambda}\to{\text{limit}}  ~\theta(\alpha)$
  exists, then so is  ~  $\underset{\alpha\in
 \Lambda}\to{\text{limit}}  ~(p\theta(\alpha))$ , for any fixed $p\in R(T)$.
Moreover
$$ \underset{\alpha\in
 \Lambda}\to{\text{limit}}  ~(p\theta(\alpha)) = p ~\underset{\alpha\in
 \Lambda}\to{\text{limit}}  ~\theta(\alpha)~. \tag"(1)"$$

 $$\text{ch}_T \Bbb C [v^{-1}\overline{BwB} \cap U^-]
 = \underset{\lambda\in D}\to{\text{limit}} ~
(\delta_{v^{-1}} \cdot (e^{v\lambda} \ast (y_w\cdot e^\lambda))).
 \tag"\text{{\bf (2.11) Corollary.}}"
 $$

 \demo{Proof}  By the previous proposition and the Demazure character formula
(cf.  \cite{A}, \cite{Jo2}, \cite{Ra2, Remarks 4.4}, \cite{Se}, \cite{Ku,
 Theorem 3.4}, \cite{Ma}),
 $$\align
 \text{ch}_T \Bbb C [v^{-1}\overline{BwB} \cap U^-]
 &=  \underset{\lambda\in D}\to{\text{limit}} ~
(\delta_{v^{-1}} \cdot (e^{v\lambda} \text{ch}_T (V_w(\lambda)^\ast )))  \\
 &= \underset{\lambda\in D}\to{\text{limit}} ~
 (\delta_{v^{-1}} \cdot (e^{v\lambda} \ast (y_w\cdot e^\lambda))).
 \endalign$$

  Observe that the existence of the above limit is
 guaranteed by Proposition (2.9) and the fact that
 $\Bbb C[v^{-1}\overline{BwB} \cap U^-]$ is an admissible
 $T$-module (being quotient of $\Bbb C[U^-]$).\qed

  \enddemo
 \vskip1ex

 Finally we come to the proof of Theorem (2.2).

 \demo{ {\bf \S(2.12)} Proof of Theorem (2.2)}

 Write (cf. (4) of \S2.1)
 $$
 y_w = \sum_{u\leq w}
 b_{w^{-1},u^{-1}}\delta_u~. $$
 Then
 $$ \align e^{v\lambda}\ast (y_w\cdot e^\lambda) &=
 \sum_u(\ast b_{w^{-1},u^{-1}}) e^{v\lambda-u\lambda}\\
 &= \ast b_{w^{-1},v^{-1}} + \sum\Sb u\neq v\\u\leq
 w\endSb
 (\ast b_{w^{-1},u^{-1}} ) e^{v\lambda-u\lambda} .
 \tag1 \endalign
 $$
 For any (regular)  weight  $\lambda_o \in D^o$ ,  $v\lambda_o-
u \lambda_o \neq 0$ for $u \neq v$.
 From the definition of $b_{w^{-1},u^{-1}}$, it is easy
 to see that there exist positive roots $\{ \beta_1,\dots
 ,\beta_\ell  \}$ depending on $w$ (possibly with
 repetitions) such that
 $$
  P \ast b_{w^{-1},u^{-1}} \in
 R(T) \qquad \text{for all } u\leq w~, \tag2
 $$
 where $P:=\prod^\ell _{k=1}
 (1-e^{-\beta_k }) $ .
\vskip1ex
Fix $\lambda_o \in D^o$. Then the subset $\{n\lambda_o\}_{n\geq 1} \subset D$
being cofinal in $D$ under $\prec$ ,
$$\underset{\lambda\in D}\to{\text{limit}} ~
  (e^{v\lambda} \ast (y_w\cdot e^\lambda))
=\underset{n\to \infty}\to{\text{limit}} ~
  (e^{nv\lambda_o} \ast (y_w\cdot e^{n\lambda_o})) ~.\tag"(3)" $$
Then by (1) of Definition (2.10) and (1), (2), (3) as above, we get

 $$\align   P~\underset{\lambda\in D}\to{\text{limit}} ~
  (e^{v\lambda} \ast (y_w\cdot e^\lambda))
&=\underset{n\to \infty}\to{\text{limit}} ~
  (P(e^{nv\lambda_o} \ast (y_w\cdot e^{n\lambda_o})))\\
&=  \underset{n\to \infty}\to{\text{limit}} ~
  (P\ast b_{w^{-1},v^{-1}} + \sum\Sb u\neq v\\u\leq
 w\endSb
 (P\ast b_{w^{-1},u^{-1}} ) e^{n(v\lambda_o-u\lambda_o)})
\\
 &=P\ast b_{w^{-1},v^{-1}} + \sum\Sb u\neq v\\u\leq
 w\endSb
 (P\ast b_{w^{-1},u^{-1}} ) ~\underset{n\to \infty}\to{\text{limit}}
{}~(e^{n(v\lambda_o-u\lambda_o)}) \\
&=P\ast b_{w^{-1},v^{-1}} ~.
  \endalign$$
So, we get (in the  $Q(T)$-module $\widetilde{Q(T)}$)
$$1\otimes \underset{\lambda\in D}\to{\text{limit}} ~
  (e^{v\lambda} \ast (y_w\cdot e^\lambda)) =\ast b_{w^{-1},v^{-1}} ~. \tag"(4)"
$$
 So, by Corollary (2.11)  and Identity (4),
 we get
 $$
 \text{ch}_T \Bbb C[v^{-1}\overline{BwB}\cap U^-]
 =\delta_{v^{-1}}\cdot( \ast b_{w^{-1},v^{-1}}).
 $$
 But the variety $v^{-1}\overline{BwB}\cap U^-$ provides  an affine
 neighborhood of the point  $\frak e\in v^{-1}X_w $.  In
 particular,
$$
 \text{gr}\, \Cal O_{\frak e,v^{-1}X_w} \cong \text{gr}\,
 {\Bbb C}[v^{-1} \overline{BwB} \cap U^-].
 $$

 The theorem now follows from the complete reducibility
 of the $T$-module ${\Bbb C}[v^{-1}\overline{BwB}\cap U^-]$, by translating the
variety $v^{-1}X_w $ under $\overline{v}.\qed $
 \enddemo

\flushpar
{\bf (2.13)} {\it Remarks.}
(1) This theorem was obtained by the author in 1987 and privately
circulated in  the  preprint `` A connection of equivariant $K$-theory
with the singularity of Schubert varieties".

(2) A different proof of the Theorem was subsequently given by Bressler [Br].
Even though I have not seen, M. Brion mentioned to me that he also obtained a
proof of this theorem (unpublished).
\vskip8ex
\centerline{{\bf 3. Some consequences of Theorem (2.2)}}

\vskip3ex
 After the following definitions, we give some of the
 corollaries of Theorem (2.2).
 \vskip2ex

 \flushpar {\bf (3.1) Definitions.} (a) For any $\ell \in
 \Bbb Z^+:= \{ 0,1,2,\dots  \}$ and any $a=\sum n_\lambda
 e^\lambda\in R(T)$, denote by $(a)_\ell =\sum n_\lambda
 \frac{\lambda^\ell }{\ell !}\in S^\ell (\frak{h}^\ast )$ , where
$S^\ell (\frak{h}^\ast )$
 is the space of homogeneous polynomials of degree $\ell
 $ on
 $\frak{h}$.  Further, denote by $[a]= (a)_{\ell _0}$
 where $\ell _0$ is the smallest element of $\Bbb Z^+$
 such that $(a)_{\ell _0}\neq 0$.  (If $a$ itself is 0,
 we define $[ a ]=0$.) Now for $q=\frac{a}b \in Q(T)$,
 where $a,b\in R(T)$, we define $[q]=\frac{[a]}{[b]}\in
 Q(\frak{h})$ (the quotient field of the symmetric
 algebra $S(\frak{h}^\ast )$).  Clearly $[ q ]$ is well
 defined.

 When $q\neq 0$ and deg $[ a ]\leq $ deg $[ b ]$, we say
 that $q$ has a pole (at the identity $e$) of order = deg
 $[ b ]-$ deg $[ a ]$.  It is easy to see that
 $b_{w^{-1},v^{-1}}$ (cf. (4) of \S2.1), when non-zero,
 has a pole of order $\leq \ell (w)$.
 \vskip1ex
 \flushpar (b) {\it The nil Hecke ring $Q_W$} (\cite{KK1,
 \S4}):  Let $Q_W$ be the smash product of the $W$-field
 $Q(\frak{h})$ with the group algebra $\Bbb Z[W]$, with
 the product given by the same formula (1) in \S2.1.  For
 any simple reflection $r_i,\ 1\leq i\leq n$, define
 $x_{r_i}\in Q_W$ by $x_{r_i}=
 -(\delta_{r_i}+\delta_e)\frac{1}{\alpha_i}$.  Now, for
 any $w\in W$, define $x_w= x_{r_{i_1}}\dots x_{r_{i_p}}$ ,
where $w=r_{i_1} \dots r_{i_p}$
 is a reduced decomposition.  The element $x_w$ is well
 defined by \cite{KK1, Proposition 4.2}.  Write, as in
 \cite{KK1, Proposition 4.3},
 $$
 x_w = \sum_v c_{w^{-1},v^{-1}} \ \delta_v,\ \text{for
 some (unique)}
 c_{w^{-1},v^{-1}} \in Q(\frak{h}).
 $$

 \proclaim{(3.2) Corollaries  (of Theorem 2.2)}  For any
 $v,w\in W$:
 \roster
 \item "{\text (a)}" $b_{w^{-1},v^{-1}}\neq 0$ if and only if $v\leq
 w$; and in this case it has a pole of order exactly
 equal to $\ell (w)$. Further,
 $$
 \biggl( \prod_{\beta \in \Delta _+} (1-e^\beta ) \biggr)
 b_{w^{-1},v^{-1}}
 \in R(T).
 \tag"(1)"
 $$
 \item "{\text (b)}" $[\ast b_{w^{-1},v^{-1}} ] =
 c_{w^{-1},v^{-1}}$; and hence for any $v\leq w$,
$$
 [\text{{\text ch}}\,(\text{{\text gr}}\,\Cal O_{\frak v,X_w})] =
 c_{w^{-1},v^{-1}} ~,$$
 as elements of $Q(\frak{h})$.

{\it In particular, $c_{w, v} \neq 0$ if and only if} $v
\leq w$.

 {\it Further}
$$({\displaystyle{\prod_{\beta \in
\triangle_{+}}}}\beta) c_{w,v} \in S ({\frak h}^*).  \tag"(2)"$$
\endroster
 \endproclaim

 \demo{Proof}  As observed in \S2.1(d), $b_{w,v}=0$ unless
 $v\leq w$.  So let us assume that $v\leq w$. Set $\Cal A^v= vU^-\frak e
\subset G/B$ .
 Since $\Cal A^v\cap X_w$ is  a closed
 subvariety of the affine space $\Cal A^v$, $\text{Tor}^p_{\Bbb C[\Cal A^v]} (
\Bbb C[\Cal A^v\cap X_w] , \Bbb C)$ is a finite dimensional vector space over
$\Bbb C$ for any $p$ and moreover ($\Cal A^v$ being smooth) is $0$ for large
enough $p$.  Set
$$ F= \sum_p (-1)^p \text{ch} (\text{Tor}_p^{\Bbb C[\Cal A^v]}
(\Bbb C[\Cal A^v\cap X_w], \Bbb C) ) \in R(T). $$
Then from the Koszul complex we get,
 $$
 (\prod_{\beta\in \Delta _+}
 (1-e^{v\beta}) ) \text{ch}\, \Bbb C[\Cal A^v\cap X_w] =
  F , ~\text{as elements of}~ R(T). \tag"(3)"$$
 It can be easily seen that the coefficient of $e^0$ in the left side of the
above identity is non-zero, in particular $F\neq 0$. From (3) we obtain $$
\align  1\otimes \text{ch}\, \Bbb C[\Cal A^v\cap X_w] &= F \cdot
\prod_{\beta\in \Delta _+}
 (1-e^{v\beta})^{-1}  \\
 &= F \prod_{\gamma \in \Delta _+\cap v\Delta_- } (-e^\gamma)\
 \prod_{\beta\in \Delta _+} (1-e^{\beta})^{-1}, \, \text{as
 elements of }
\widetilde{Q(T)} . \tag"(4)"
 \endalign
 $$

 From (4) it is clear that $1\otimes $ ch$\,\Bbb C
 [\Cal A^v\cap X_w]\neq 0$ as an element of $\widetilde{Q(T)}$.
Moreover, since  $\Cal A^v\cap X_w$ is an affine neighborhood of $\frak v$ in
$X_w$, we get
$$  \text{ch}\, \Bbb C[\Cal A^v\cap X_w] = \text{{\text ch}}\,(\text{{\text
gr}}\,\Cal O_{\frak v,X_w}) . \tag"(5)" $$
 But then by (5) and Theorem (2.2), we get that
 $b_{w^{-1},v^{-1}} \neq 0$.  The assertion that
 $b_{w^{-1},v^{-1}} $ has a pole of order exactly equal
 to $\ell (w)$ (whenever $v\leq w$) follows from a lemma
 of Joseph \cite{Jo, \S2.3}.  This
 proves the first part of Corollary (a).  Assertion (1)
 of  part (a) follows immediately from (4) (and
 Theorem 2.2).

 To prove part (b), in view of Theorem (2.2), we only need
 to show that
 $$
 [ \ast b_{w^{-1},v^{-1}} ] =
 c_{w^{-1},v^{-1}} :
 \tag"(6)"
 $$

 By induction (on $\ell (w)$) we assume the validity of
 (6) for any $w$ with $\ell (w)\leq k$ and any
 $v\in W$ and take $w'=wr_i$ of length $k+1$, where $r_i$
 is a simple reflection such that $\ell (w')>\ell (w)$.
 (The case $w=e$ is obviously true.)

  By Definition 2.1 (d),$$\align
 y_{w'}&:=
  y_wy_{r_i} \\
 &=  (\sum_v b_{w^{-1},v^{-1}} \delta_v)
(\delta_e+\delta_{r_i})(\frac{1}{1-e^{-\alpha_i}})\\
 &= \sum_v \frac{ b_{w^{-1},v^{-1}}+b_{w^{-1},r_iv^{-1}}}
{1-e^{-v\alpha_i}} \delta_v.
 \endalign
 $$
This gives for any $v\in W$,
$$b_{r_iw^{-1},v^{-1}}=\frac{b_{w^{-1},v^{-1}}+b_{w^{-1},r_iv^{-1}}}
{1-e^{-v\alpha_i}} . \tag"(7)"
 $$

  Exactly the
 same way, using the definitions from \S3.1(b), we
 obtain:
 $$
 c_{r_iw^{-1},v^{-1}}
 =
 \frac{c_{w^{-1},v^{-1}} +c_{w^{-1},r_iv^{-1}}
 }{-v\alpha_i} .
 \tag"(8)"
 $$

 By (7) and part (a) of the corollary we get:
 $$
 [ \ast b_{r_iw^{-1},v^{-1}} ]
 =
 \frac{ [ \ast b_{w^{-1},v^{-1}} ] + [
 \ast b_{w^{-1},r_iv^{-1}} ] }{-v\alpha_i}.
 $$

 Hence by the induction hypothesis (using  (8)),
 (6) follows for $w'=wr_i$.  This completes the
 proof of Corollaries (3.2).\qed
 \enddemo

 \flushpar {\bf (3.3) Remarks.}  (1)  The (b)-part of the
 above corollary is due to Rossmann \cite{R, \S3.2}.  In
 fact, this motivated our theorem (2.2).
 \par (2) The assertions (1) and (2) as above can be derived purely
 algebraically (cf. \cite{KK2, Corollary 4.18 and Remark 4.17(b)
 }).
\vskip2ex

 The following lemma gives an expression for $b_{w,v}$
 (and $c_{w,v}$) and can be easily proved by using the
 definitions.

 \proclaim{(3.4) Lemma} Fix any $v\leq w\in W$, and take
 a reduced decomposition $w=r_{i_1}\dots r_{i_p}$.
 Then
 $$\gather
 b_{w^{-1},v^{-1}} =
  \sum((1-e^{-r^{\epsilon_1}_{i_1}\alpha_{i_1}})
 ( 1-e^{ -r^{\epsilon_1}_{i_1} r^{\epsilon_2}_{i_2}\alpha_{i_2}} )
 \cdots
 (1-e^{ -r^{\epsilon_1}_{i_1}
 \dots r^{\epsilon_p}_{i_p}\alpha_{i_p}}))^{-1}.\\
 \intertext{Similarly}
 c_{w^{-1},v^{-1}} = (-1)^p \sum
 ((r^{\epsilon_1}_{i_1}\alpha_{i_1})
 ( r^{\epsilon_1}_{i_1} r^{\epsilon_2}_{i_2}\alpha_{i_2} )
 \cdots
 ( r^{\epsilon_1}_{i_1} \dots r^{\epsilon_p}_{i_p}\alpha_{i_p}))^{-1},
 \endgather$$
 where both the sums run over all those $(\epsilon_1,\dots
 ,\epsilon_p)\in \{ 0,1\}^p$ satisfying $r^{\epsilon_1}_{i_1} \dots
 r^{\epsilon_p}_{i_p}=v$.  (The notation $r^0_i$ means the
 identity element.)
 \endproclaim

 \vskip7mm
\newpage
 \centerline{{\bf 4. Ring of functions on the tangent
 cone --}} \centerline{{\bf the graded algebra structure}}
 \vskip1ex

 \flushpar
{\bf \S (4.1)}  For any $\lambda\in D$, the (finite
 dimensional) $G$-module $V(\lambda)$ admits a filtration
 $\{ \Cal F_p(\lambda) \}_{p\geq 0}$ as follows:
 \par Let $\{ U_p(\frak{u}^-\}_{p\geq 0}$ be the standard
 filtration of the universal enveloping algebra
 $U(\frak{u}^-)$, where we recall that $U_p(\frak{u}^-)$
 is the span of the monomials $X_1\dots X_m$ for
 $X_i\in \frak{u}^-$ and $m\leq p$.  Now set
 $$
 \Cal F_p(\lambda)= U_p(\frak{u}^-)\cdot e_\lambda,
 $$
 where $e_\lambda$ is any non-zero highest weight vector
 in $V(\lambda)$.

 Fix $\lambda\in D$, $v\leq w\in W$ ,  $\theta\in
 V(\lambda)^\ast $, and a highest weight vector
 $e_\lambda\in V(\lambda)$.  Recall the definition of  the function
$\varphi_\lambda$ from \S2.4. we abbreviate $\varphi_\lambda (\theta \otimes
e_\lambda)$ by $\varphi^\theta$.  Thus   $\varphi^\theta :  U^-\to \Bbb C$ is
the function
 $$
 \varphi^\theta(g) = \langle \theta,ge_\lambda \rangle,\quad
 \text{for }g\in U^-. $$
 By Lemma (2.4), $\varphi^\theta$ vanishes on
 $v^{-1}\overline{BwB}\cap U^-\Leftrightarrow \langle
 \theta,v^{-1}V_w(\lambda) \rangle=0$.  Identify $U^-$
 with the affine space $\frak{u}^-$ under the exponential
 map.  This gives rise to a gradation on $\Bbb C[U^-]$.
 Now let $\varphi^\theta_{d}$ be the $d^{\text{th}}$
 graded component of $\varphi^\theta$ (for any $d\geq 0$),
 i.e.,
 $$
 \varphi^\theta_d(X) =
 \frac{1}{d !} \langle \theta,X^d \cdot
 e_\lambda\rangle,\quad \text{for }X\in \frak u^-.
 \tag"$(\ast)$"
 $$

 The following  lemma follows immediately from $( \ast
 )$, if we use the fact that for any vector space $V$,
 its $p^{th}$ symmetric power $S^p(V)$ is spanned by $\{
 v^p \}_{v\in V}$.

 \proclaim{(4.2) Lemma} Fix $p\geq 1$.  Then for any $\theta$
 as above (i.e. $\langle \theta,v^{-1}V_w(\lambda)
 \rangle=0$), $\varphi^\theta_d\equiv 0$ for all $0\leq
 d<p$ if and only if $\langle
 \theta,v^{-1}V_w(\lambda) + \Cal
 F_{p-1}(\lambda)\rangle=0$.\qed
 \endproclaim

 For any $p\geq 1$ and any subvariety $0\in Y\subset \Bbb
 A^n$, let $\Cal I_p(Y)$ denote the set of degree
 $p^{\text{th}}$ components of all those functions $f$ in
 the ideal $\Cal I(Y)$ of $Y\subset \Bbb A^n$, such that
 the $d^{\text{th}}$ homogeneous component $f_d$ of $f$
 is 0 for all $d<p$.
 \vskip1ex

 As an immediate consequence of the above lemma, we get
 the following
 \proclaim{(4.3) Corollary} For any $p\geq 0$ the map
 $
 \theta\otimes e_\lambda\mapsto  (\varphi_\lambda (\theta \otimes
e_\lambda))_p$
 induces a $T$-equivariant  injective map
 $$
 \biggl(
 \frac{v^{-1}V_w(\lambda)+ \Cal
 F_p(\lambda)}{v^{-1}V_w(\lambda)+\Cal F_{p-1}(\lambda)}\biggr)^\ast
 \otimes  \Bbb C_\lambda \overset
 f_p(v,w;\lambda)\to{\hookrightarrow}
 \Cal I_p (v^{-1}\overline{BwB} \cap U^-) ,
 $$
 where $\Bbb C_\lambda\subset V(\lambda)$ is the highest
 weight subspace, and $\Cal F_{-1}(\lambda)$ is defined
 to be 0.\qed
 \endproclaim

 It is easy to see that under the map
 $
 \hat{\delta}_{\lambda,\mu } : V(\lambda+\mu )
 \longrightarrow
 V(\lambda)\otimes \Bbb C_\mu $(of \S 2.6), the image
 $\hat{\delta}_{\lambda,\mu } (\Cal F_p(\lambda+\mu ))
 \subset \Cal F_p(\lambda) \otimes \Bbb C_\mu .$
 Moreover, by Lemma (2.7),
$$ \hat{\delta}_{\lambda,\mu } (v^{-1}V_w(\lambda+\mu ))
 \subset v^{-1}V_w (\lambda)\otimes \Bbb C_\mu .
 $$

 In particular, $\hat{\delta}_{\lambda,\mu }$ gives rise to a
 $T$-module map \, $\delta_{\lambda,\mu }(v,w;p)$ making
 the following diagram commutative (for any $\lambda,\mu
 \in D,\ v\leq w\in W$ and $p\geq 0$):
  $$
   \gather
  \biggl( \frac{v^{-1} V_w(\lambda)+\Cal
 F_p(\lambda)}{v^{-1}V_w(\lambda)+\Cal F_{p-1}(\lambda)}
 \bigg)^\ast \otimes \Bbb C_\lambda
  \quad
   \overset{\delta_{\lambda,\mu
 }(v,w;p)}\to\longrightarrow
  \quad
  \biggl( \frac{v^{-1} V_w(\lambda+\mu)+\Cal
 F_p(\lambda+\mu)}{v^{-1}V_w(\lambda+\mu )+\Cal
 F_{p-1}(\lambda+\mu )} \bigg)^\ast \otimes \Bbb
 C_{\lambda+\mu} \\
 \vspace{4\jot}
   \sideset{f_p(v,w;\lambda)}\and\to\searrow
  \qquad\qquad\qquad
  \sideset\and{f_p(v,w;\lambda+\mu)}\to\swarrow\\
 \vspace{2\jot}
   \Cal I_p (v^{-1} \overline{BwB} \cap U^-) .
   \endgather
  $$
 By the injectivity of the map $f_p(v,w;\lambda)$, we see
 that $\delta_{\lambda,\mu }(v,w;p)$ is injective.  Thus
 $$
 \{ \biggl( \frac{v^{-1} V_w(\lambda)+\Cal
 F_p(\lambda)}{v^{-1}V_w(\lambda)+\Cal F_{p-1}(\lambda)}
 \bigg)^\ast \otimes \Bbb C_\lambda \biggr\}_{\lambda\in
 D}
 $$
 forms a directed system of $T$-modules and there is an
 induced $T$-module map
 $$
 f_p(v,w): \underset{\lambda\in D}\to{\text{limit}_\to}
 \biggl(
 \biggl[
 \frac{v^{-1} V_w(\lambda)+\Cal
 F_p(\lambda)}{v^{-1}V_w(\lambda)+\Cal F_{p-1}(\lambda)}
 \bigg]^\ast \otimes \Bbb C_\lambda
 \biggr)
 \rightarrow
 \Cal I_p (v^{-1}\overline{BwB} \cap U^-).
 $$

 \proclaim{(4.4) Theorem} The above map $f_p(v,w)$ is a
 $T$-equivariant isomorphism for all $p\geq 0$ and all
 $v\leq w\in W$.  In particular, there is a
 $T$-equivariant isomorphism
 $$
 \text{{\text gr}}_p (\Cal O_{\frak{e},v^{-1}X_w})
 \approx \underset{\lambda\in D}\to{\text{limit}_\to}
 \biggl(  \biggl[ \frac{v^{-1} V_w(\lambda)\cap
 \Cal F_p(\lambda)}{v^{-1}V_w(\lambda)\cap \Cal F_{p-1}(\lambda)}
 \bigg]^\ast \otimes \Bbb C_\lambda  \biggr),
 \tag1
 $$
 where {\text gr}$_p(\Cal O_{e,v^{-1}X_w})$ is the
 p$^{\text{th}}$ graded component of {\text
 gr}$(\Cal O_{\frak e,v^{-1}X_w}).$
 \endproclaim

 \demo{Proof} Since $f_p(v,w;\lambda)$ is injective for
 all $\lambda\in D$, $f_p(v,w)$ is clearly injective.
 The surjectivity of $f_p(v,w)$ follows from Lemma (2.4)
 and Lemma (4.2).
 \vskip1ex

 We now come to the proof of (1):  Observe first that by
 \cite{Ha, Lecture 20},
 $$
 \text{{\text gr}}_p (\Cal O_{\frak e ,v^{-1}X_w})  \approx
 S^p((\frak u^-)^\ast )/ \Cal
 I_p(v^{-1}\overline{BwB}\cap U^-),
 \tag2$$
 where $S^p$ is the p$^{\text{th}}$ symmetric power.

 Now for any (fixed) $p$, if we take $\lambda$ to be
 sufficiently large, then the map
 $$
 U_p (\frak{u}^-)
  \otimes \Bbb C_\lambda \to
 \Cal F_p(\lambda) \text{ given by }X\otimes e_\lambda
 \mapsto Xe_\lambda,
 $$
 for $e_\lambda\in \Bbb C_\lambda\subset V(\lambda)$, is
 a $T$-module isomorphism.  In particular, by the
 Poincare-Birkhoff-Witt theorem,
 $$
 \Cal F_p(\lambda) / \Cal F_{p-1}(\lambda)
 \approx S^p(\frak{u}^-)\otimes \Bbb C_\lambda
 \quad\text{(for large enough $\lambda$).}
 \tag3
 $$

 Consider the exact sequence
 $$
 0
 \rightarrow
 \frac{v^{-1} V_w (\lambda)\cap \Cal F_p(\lambda)}{v^{-1}
 V_w (\lambda) \cap \Cal F_{p-1}(\lambda)}
 \rightarrow
 \frac{\Cal F_p(\lambda)}{\Cal F_{p-1}(\lambda)}
 \rightarrow
 \frac{v^{-1} V_w (\lambda) + \Cal F_p (\lambda)}{v^{-1}
 V_w (\lambda) + \Cal F_{p-1}(\lambda)}
 \rightarrow 0.
 $$
 Dualizing this sequence and using (3) we get (for large
 enough $\lambda$)
 $$
 0
 \leftarrow
 \biggl( \frac{v^{-1} V_w (\lambda)\cap \Cal
 F_p(\lambda)}{v^{-1} V_w (\lambda) \cap \Cal
 F_{p-1}(\lambda)}\biggr)^\ast \otimes \Bbb C_\lambda
 \leftarrow S^p((\frak{u}^-)^\ast) \leftarrow
 \biggl( \frac{v^{-1} V_w (\lambda) + \Cal F_p
 (\lambda)}{v^{-1} V_w (\lambda) + \Cal F_{p-1}(\lambda)}\biggr)^\ast
\otimes {\Bbb C}_\lambda
 \leftarrow 0.
 $$

 Now the isomorphism (1) is established from (2) and the
 isomorphism $f_p(v,w)$. \qed
 \enddemo
\flushpar
 {\bf \S(4.5)} For any variety $X$ and a closed point $x\in X$, let
 $Z_x(X)$ denote the Zariski tangent space of $X$ at $x$.
 For any subvariety $Y\subset G/B$ containing the base
 point $\frak{e}$, we get the induced inclusion
 $Z_{\frak{e}}(Y)\hookrightarrow Z_{\frak{e}}(G/B)$.  But
 $Z_{\frak{e}}(G/B)$ can be canonically identified with
 $\frak{u}^-$ (since $U^-$ is an open neighborhood
 around $\frak{e}$ in $G/B$), in particular,
 $Z_{\frak{e}}(Y)$ can be canonically viewed as a
 subspace of $\frak{u}^-$.  For any $\alpha\in \Delta $,
 let $r_\alpha\in W$ be the reflection defined by
 $r_\alpha(\lambda)= \lambda- \langle \lambda,\alpha^\vee
 \rangle \alpha$.
 \vskip1ex

The following result is due to Polo \cite{P, Theorem 3.2}. It may be recalled
that a  different description of the   Zariski tangent
 space in the case of classical groups was given by Lakshmibai-Seshadri (cf.
\cite{LS} \cite{L}).  Observe that by virtue of the
 automorphism of $G/B$, given by $gB\mapsto \overline{v}
 gB$ (for  $g\in G$), $Z_{\frak{e}}(v^{-1}X_w)$ is
 isomorphic with $Z_{\frak{v}}(X_w)$.

The first part of the following result follows immediately from Theorem (4.4)
and the second part follows from the fact that $X_w \subset
G/B$  is defined by linear equations.

 \proclaim{(4.6) Corollary} For any $v\leq w$,
 $$\gather
 Z_{\frak{e}}(v^{-1}X_w) = \{ X\in \frak{u}^- :
 Xe_\lambda\in v^{-1}V_w(\lambda), \ \text{for all }
 \lambda\in D\}.
 \tag1\\
 \intertext{In fact,}
 Z_{\frak{e}}(v^{-1}X_w) = \{ X\in \frak{u}^- :
 Xe_{\lambda_0} \in v^{-1} V_w (\lambda_0)\}, \ \text{for
 any one regular  } \lambda_0\in D,
 \tag2
 \endgather
 $$
where $e_\lambda$ is a non-zero highest weight vector of $V(\lambda)$.
 \endproclaim

 \demo{Proof}  The identity (1) follows from Theorem (4.4)
 immediately, since $Z_{\frak e}(v^{-1}X_w)= $
 $\text{gr}_1(\Cal O_{\frak{e},v^{-1}X_w})^\ast $.  However, we give
 the following  direct proof:

 Fix $\lambda\in D$ and take $\theta\in
 (V(\lambda)/ v^{-1}V_w(\lambda))^\ast$, and
 consider the corresponding function $\varphi^\theta: U^-\to \Bbb C$
 defined by
 $$
 \varphi^\theta(\exp X)= \langle \theta,\exp Xe_\lambda \rangle =
 \langle \theta,Xe_\lambda \rangle + \text{ order
 two and higher terms}.
 $$
(Observe that  $\langle \theta,e_\lambda
 \rangle=0$ by assumption.)
 So the linear part $L(\varphi^\theta)\in (\frak{u}^-)^\ast $ (under
 the identification exp :$\frak{u}^-\to U^-$) of $\varphi^\theta$
 is given by
 $$
 L(\varphi^\theta)X= \langle \theta,Xe_\lambda  \rangle, \quad
 \text{for }X\in \frak{u}^-.
 \tag3$$

 Let $\Cal I(v^{-1}\overline{BwB}\cap U^-)$ denote the ideal of the
 closed subvariety $v^{-1}\overline{BwB}\cap U^-$ of $U^-$.  Then, by
 the definition of the Zariski tangent space,
 $$\spreadlines{2\jot}\align
 Z_{\frak{e}} (v^{-1}X_w)
 &=  \{ X\in \frak{u}^- : L(f)X=0, \text{ for all }f\in
 \Cal I(v^{-1}\overline{BwB}\cap U^-)\},\\
 &=  \{ X\in \frak{u}^- :Xe_\lambda\in
 v^{-1}V_w(\lambda), \text{ for all }\lambda\in D\},
 \text{ by (3) and Lemma (2.4)}.
 \endalign$$
 This proves (1).

We now prove (2): The tensor product of  sections gives rise to
an algebra structure on the space $R:= \oplus_{m\geq 0} H^0(G/B, \frak
L(m\lambda_o))$. Let $K_m$ be the kernel of the restriction map
 $ H^0(G/B, \frak L(m\lambda_o)) \to H^0(X_w, \frak L(m\lambda_o)\vert_{X_w})$.
Then by a result of Ramanathan \cite{Ra2,
Theorem 3.11}, the kernel $K :=\sum_{m\geq 0} K_m$  of the surjective map  $
\oplus_{m\geq 0}  H^0(G/B, \frak L(m\lambda_o)) \to \oplus_{m\geq 0}  H^0(X_w,
\frak L(m\lambda_o)\vert_{X_w})$ is generated as an ideal in the ring $R$ ) by
$K_1$ (i.e. $X_w$ is linearly defined in $G/B$ with respect to $\frak
L(\lambda_o)$). This, in particular,
implies (by translating via $\overline{v}^{-1}$ and using Lemma 2.3) that the
ideal  $\Cal I(v^{-1}\overline{BwB}\cap U^-)$ is generated
by the functions $\{\varphi^\theta\}$ where $\theta$ ranges over
$ (V(\lambda_o)/ v^{-1}V_w(\lambda_o))^\ast$.  Now by an argument identical to
the proof of (1), we get (2). \qed
\enddemo

 \proclaim{(4.7) Lemma}  Let $\frak{g}$ be simply-laced.
 Assume that there exist integers $p,\ p_1,\dots,p_k\geq
 1$ and roots $\beta ,\beta _1,\dots ,\beta _k\in \Delta
 _+$ such that
 $$\gather
 p\beta = \sum^k_{j=1} p_j\beta j
 \tag1\\
 \intertext{and}
 \sum p_j\leq p.\tag2
 \endgather
 $$
 Then $\beta j= \beta  $ ,  for all $1\leq j\leq k$.
 \endproclaim

 \demo{Proof} We can assume without loss of generality
 that no $\beta j= \beta$.  Now by (1) we get
 $$
 p\langle \beta ,\beta ^\vee \rangle =
 \sum^k_{j=1} p_j
 \langle \beta j,\beta ^\vee \rangle .
 \tag3
 $$
 But $\frak g$ being simply-laced, $\langle \beta j,\beta ^\vee
 \rangle\leq 1$ (since $\beta j\neq \beta $), and hence
 by (2) and (3) we get
 $$
 2p \leq
 \sum^k_{j=1} p_j  \leq p.
 $$
 This contradiction proves the lemma.\qed
 \enddemo

 \flushpar {\bf (4.8) Definition.} For any $v \leq w \in W$, define  $S(w,v)=
\{ \alpha\in \Delta _+: vr_{\alpha}
 \leq w\}$.
Then as is easy to see $\#S(w,v) = \#S(w^{-1},v^{-1})$.

 \proclaim{(4.9) Proposition} Let $\frak{g}$ be simply-laced.
 Fix $v\leq w\in W$.  Then for any $\alpha\in \Delta _+$
 such that $\alpha\notin S(w,v)$ but $E_{-\alpha}\in
 Z_{\frak e}(v^{-1}X_w)$, there exists a non-zero element
 $\theta_\alpha\in  \text{gr}_1 (\Cal O_{\frak{e},v^{-1}X_w})$ of weight
$\alpha$
 satisfying $\theta_\alpha^{\langle \rho,\alpha^\vee
 \rangle}=0$
  as an element of {\text gr}$_{\langle \rho,\alpha^\vee
 \rangle} (\Cal O_{\frak{e},v^{-1}X_w})$, where $E_{-\alpha}$
 is a non-zero root vector of $\frak g$ corresponding to the negative root
$-\alpha$, and $\rho$ is the half sum of positive roots.
 In particular, the tangent cone
 $T_{\frak{e}}(v^{-1}X_w)$ is non-reduced in this case.
 \endproclaim

 \demo{Proof} By Lemma (4.7), the weight space of
 $U_p(\frak{u}^-)$ corresponding to the weight  $-p\alpha$
 (for any $p\geq 1$ ) is one dimensional, and is spanned by
 $E_{-\alpha}^p$.  Since gr$_1(\Cal O_{y,Y})$ is canonically
 isomorphic with the dual space $Z_{y,Y}^\ast $ (for any
 variety $Y$ and $y\in Y$), and $E_{-\alpha}\in
 Z_{\frak{e},v^{-1}X_w}$, there exist a non-zero element
 $\theta_\alpha\in$ gr$_1(\Cal O_{\frak{e}, v^{-1}X_w})$
 of weight $\alpha$.  Under the embedding $Z_{\frak{e},
 v^{-1}X_w}\hookrightarrow  \frak{u}^-$ (cf. \S4.5), we can identify
 the element $\theta_\alpha$ with the element of $(
 \frak{u}^- )^\ast $ defined by $\theta_\alpha(E_{-\beta
 })= \delta_{\alpha,\beta }$, for all $\beta \in \Delta
 _+$.

 By virtue of Theorem (4.4), to prove that
 $\theta^p_\alpha=0$ (where $p:= \langle \rho,\alpha^\vee
 \rangle$), it suffices to show that (for all large
 enough $\lambda\in D$)
 $$
 \theta^p_{\alpha _{\vert_{(v^{-1}V_w(\lambda)\cap
 \Cal F_p(\lambda))}}}\equiv 0 :
 $$
 Since $\theta^p_\alpha$ is of weight $p\alpha$ and the
 weight space of $U_p(\frak{u}^-)$ corresponding to the
 weight $-p\alpha$ is spanned by $E^p_{-\alpha}$, it
 suffices to show that $E^p_{-\alpha}e_\lambda\notin
 v^{-1}V_w(\lambda)$ (for all large enough $\lambda\in D$):

 For otherwise, assume that $E_{-\alpha}^p e_{\lambda_0}\in
 v^{-1}V_w(\lambda_0)$ (for some $\rho\preceq
 \lambda_0$).  Then by Lemma (2.7), $E_{-\alpha}^pe_{\rho}
 \in v^{-1}V_w(\rho)$.  But since $\bar{r }_{\alpha}
 e_\rho= E_{-\alpha}^p e_\rho$ (up to a non-zero scalar
 multiple), $\bar{r }_\alpha e_\rho \in
 v^{-1}V_w(\rho)$ and hence by \cite{BGG, Theorem 2.9}
 $vr _\alpha\leq w$, which contradicts the
 assumption and proves the proposition.\qed \enddemo

 \flushpar {\bf (4.10)}{\it  Remark.}  The `in particular' statement
 of the above Proposition can also be deduced from
 \cite{C, Theorem G(2)}.
 \vskip1ex

 For a closed  point $x$ of a scheme $X$, recall the
 definition of the tangent cone $T_x(X)$ as Spec
 (gr $\,\Cal O_x$) from \S2.1.  Define the {\sl reduced
 tangent cone} $T_x^{red}(X)$ as Spec
 (gr$^{\text{red}}\Cal O_x$), where
 gr$^{\text{red}}\Cal O_x$ =
$ (\text{gr} \,\Cal O_x )/ N$ and $N$ is the ideal
 consisting of all the nilpotent elements in
 gr $\Cal O_x$.

 The following result is due to Carrell-Peterson
 \cite{C, Theorem EG}, proved by different methods.

 \proclaim{(4.11) Corollary}  Let $\frak{g}$ be an arbitrary
 semisimple Lie algebra and fix $v\leq w\in W$.  Assume
 that $T_{\frak{e}}^{\text{red}}(\theta^{-1}X_w)$ is an
 affine space for all $v\leq \theta\leq w$.  Then the
 point $\frak{v}\in X_w$ is rationally smooth.

 Conversely, in the case when $\frak{g}$ is simply-laced,
 if the point $\frak{v}\in X_w$ is rationally smooth,
 then $T_{\frak{e}}^{\text{red}} (\theta^{-1}X_w)$ is an
 affine space for all $v\leq \theta\leq w$.
 \endproclaim
 \demo{Proof} As follows from \cite{C, Theorem F} (cf. also \cite{P,
Proposition 4.2}), for any $\alpha\in S(w,\theta)$,  $E_{-\alpha}\in
 Z_{\frak{e}}(\theta^{-1}X_w)\cong \,\text{gr}_1
 (\Cal O_{\frak{e},\theta^{-1}X_w})^\ast $.  Choose a
 non-zero element $\theta_\alpha$ of weight $\alpha$ in
 $\text{gr}_1 (\Cal O_{\frak{e},\theta^{-1}X_w})$.  Then
 $\theta ^p_\alpha\neq 0$ in $\text{gr}_p
 (\Cal O_{\frak{e},\theta^{-1}X_w})$ (for any $p\geq 1$):
 To prove this, it suffices to show that $E_{-\alpha}^p
 e_\lambda \in \theta^{-1}V_w(\lambda)$, for any $\lambda\in D$
 such that $p\leq \langle \lambda,\alpha^\vee \rangle$
 (cf. proof of Proposition 4.9):

 By the $sl(2)$-theory, $E_{-\alpha}^{\langle
 \lambda,\alpha^\vee \rangle} e_\lambda = \bar{r
 }_{\alpha}e_\lambda$ (up to non-zero scalar multiples).
If $\theta\alpha \in \Delta_-$, clearly
 $E_{-\alpha}^p e_\lambda\in \theta^{-1}V_w(\lambda)$. So assume that
$\theta\alpha \in \Delta_+$. Then (upto non-zero scalar multiples)
$$E_{-\alpha}^p e_\lambda = E_{\alpha}^{\langle
 \lambda,\alpha^\vee \rangle -p}   E_{-\alpha}^{\langle
 \lambda,\alpha^\vee \rangle } e_\lambda = E_{\alpha}^{\langle
 \lambda,\alpha^\vee \rangle -p} \bar{r}_\alpha e_\lambda \in
\theta^{-1}V_w(\lambda)~,$$
 thereby proving the claim.

 We come to the proof of the first part of the Corollary.
 Since the dimension of the tangent cone is the same as the local dimension of
the variety at that point (cf. \cite{Ha, Lecture 20}), and (by assumption)
 $T^{\text{red}}_{\frak{e}}(\theta^{-1}X_w)$ is an affine
 space, $\dim T^{\text{red}}_{\frak{e}}(\theta^{-1}X_w) =
 \ell (w)\geq \# S(w,\theta)$.  But, by Deodhar's
 conjecture (see Theorem 5.1), $\ell (w)\leq \#
 S(w,\theta)$.  Hence $\ell (w)=\#S(w,\theta)$, for all
 $v\leq \theta\leq w$.  So the first part of the
 Corollary follows from \cite{C, Theorem E}. (Observe that for any $\theta \in
W, \#\{ \alpha \in \Delta_+ : r_\alpha \theta < \theta \}=\ell (\theta)$.)

 In the simply-laced case, by Proposition (4.9) and the
 above argument,
 $$
 \dim  (\text{gr}_1^{\text{red}}
 (\Cal O_{\frak{e},\theta^{-1}X_w}))
 = \#S(w,\theta) =\ell (w),
 \tag1
 $$
 since $\frak{v}\in X_w$ is assumed to be rationally
 smooth.  But since gr$^{\text{red}}$ is generated (as an
 algebra ) by $\text{gr}_1^{\text{red}}$, we get a
 surjective map $\gamma: S( \text{gr}_1^{\text{red}}
 (\Cal O_{\frak{e},\theta^{-1}X_w}) ) \twoheadrightarrow
 \text{gr}^{\text{red}}
 (\Cal O_{\frak{e},\theta^{-1}X_w})$ (where $S$ is the
 symmetric algebra).
 But since $T_{\frak{e}}^{\text{red}}(\theta^{-1}X_w)$ is
 of dim $\ell (w)$, surjectivity of $\gamma $ and (1)
 force $\gamma$ to be an isomorphism.  This proves the
 corollary.  \qed \enddemo

 \flushpar {\bf (4.12)} {\it  Remark.}  The converse statement of the
 above corollary is not true in general for non simply- laced
 $\frak{g}$.  Take, e.g., $\frak{g}$ to be of type
 $C_2$ or $G_2$ and $w=r _1r _2r_1, ~
 v=e$. Since $\frak{g}$ is of rank 2,  (as is well known; and  can also be
proved by using Lemma 6.2 and Theorem 5.5 (a)) $\frak{e}\in X_w$
 is rationally smooth. But it can be easily seen that
 $T_{\frak{e}}^{\text{red}}(X_w)$ is not an affine space.

\vskip7ex
 \centerline{{\bf 5.  Smoothness criterion of Schubert varieties}}

 \vskip2ex

 For any $v\leq w\in W$, recall the definition of
 $S(w,v) $ from Definition (4.8).
 We recall the following very interesting conjecture of
 Deodhar \cite{D}, which was proved by Carrell-Peterson \cite{C},
 Dyer \cite{Dy}, and Polo \cite{P}.

 \proclaim{(5.1) Theorem}  For any $v \leq w\in W,\ \#
 S(w,v)\geq \ell (w)$.
 \endproclaim

 Even though the following proposition follows
 immediately by combining our Corollary 3.2(b) with
 \cite{Dy, Proposition \S3}, we give a different
 (geometric) proof (as that proof is crucially used in
 the proof of Theorem 5.5(b)).

 \proclaim{(5.2) Proposition}  Let $v\leq w\in W$.  Then
 $$
 \# S(w^{-1},v^{-1})  \,= \ell (w)
 \Leftrightarrow [ \text{ch} \, ( \text{gr} \, \Cal
 O_{\frak{v},X_w})]
 = d(-1)^{\ell (w)-\ell (v)}
 \prod_{\beta\in S(w^{-1},v^{-1})} \ \beta ^{-1},
 $$
 for some $d\in \Bbb C$.
 \endproclaim

 \demo{Proof} By Corollary 3.2 (b),
 $[\text{ch}\,(\text{gr}\, \Cal O_{\frak{v},X_w})] =
 c_{w^{-1},v^{-1}}\neq 0$ and, moreover, it can be easily
 seen from the definition of $c_{w^{-1},v^{-1}}$ that deg
 $c_{w^{-1},v^{-1}}= -\ell (w)$, where deg$\frac{P}{Q}:=
 \deg P-\deg Q$ for non-zero $P,Q\in S(\frak h^\ast) $.
 Hence the implication `$\Leftarrow$' of the above
 proposition follows.
 \vskip1ex

 Now we come to the implication `$\Rightarrow$':

 Let $\exp :\frak{u}^-\tilde\to U^- $ be the exponential map,
 where $\frak{u}^-$ is the Lie algebra of $U^-$.
 (Observe that $U^-$ being a unipotent group, $\exp$ is
 an algebraic morphism.)  Let $Y:= \exp^{-1}(U^-\frak{e}
 \cap v^{-1}X_w)$ be the closed irreducible subvariety of
 $\frak{u}^-$, where we identify $U^-$ with $U^-\frak e$.  Fix non-zero root
vectors $E_{-\beta }$
 (corresponding to the negative root $-\beta $) for
 $\beta \in \Delta _+$.  For any $\alpha \in  \Delta _+$,
 let $f_\alpha: \frak{u}^-\to \Bbb C$ be the linear map
 defined by $\sum_{\beta \in \Delta _+}t_\beta E_{-\beta
 }\mapsto t_\alpha$, and let $f^Y_\alpha$ be the
 restriction of $f_\alpha$ to $Y$.  Define the subvariety (with the
reduced structure)
 $$
 Z_S = \{ x\in Y: f^Y_\alpha(x)= 0,\ \text{for all
 }\alpha\in  S:= S(w,v)\}.
 $$
 Clearly $0\in Z_S$.  We claim that any irreducible
 component $Z_S^o$ of $Z_S$ through 0 is 0-dimensional:

 The varieties $Z_S^o \subset Z_S$ are clearly
 $T$-stable under the adjoint action of the maximal
 torus $T$ on $\frak{u}^-$.  Further, $Z_S^o$
 does not contain any 1-dimensional $T$-stable closed
 irreducible subvariety $R$:  It is easy to see that any 1-dimensional
 $T$-stable closed irreducible subvariety of $\frak{u}^-$
 is of the form $\Bbb CE_{-\beta
 }\subset \frak{u}^-$, for some $\beta \in \Delta _+$.
 In particular, $R=\Bbb C E_{-\beta_0 }$ (for some $\beta
 _0\in \Delta _+$).  This gives that $\exp (\Bbb C
 E_{-v\beta_0 }) \frak{v} \subset X_w$.  Now if $-v\beta
 _0\in \Delta _+$, then by \cite{BGG, Corollary 2.3}
 $vr _{\beta _0} < v\leq w$, so $\beta _0\in S$.  If
 $v\beta _0\in \Delta _+$, then clearly $\exp(\Bbb C
 E_{v\beta _0}) \exp(\Bbb C
 E_{-v\beta _0}) \frak{v} \subset X_w $.  In particular,
 for the subgroup $\Cal S_{vr_{\beta
 _0}v^{-1}}\subset G$ generated by $\exp (\Bbb
 CE_{-v\beta_0 })$ and $\exp (\Bbb CE_{v\beta_0 })$,
 $
 \Cal S_{vr_{\beta _0}v^{-1}} \frak{v} \subset X_w.
 $
 Again this gives,  by \cite{C, Theorem F(2)}, that $\beta
 _0\in S$.  So, in either case, $R=\Bbb CE_{-\beta_0 }$,
 for some $\beta _0\in S$.  But, by the definition of
 $Z_S$, such a $R$ is not contained in $Z_S$.  This
 contradiction establishes the claim that $Z^o_S$ does
 not contain any 1-dimensional $T$-stable closed
 irreducible subvariety.

 Embed $i:\frak{u}^-\hookrightarrow G/B$ via the map
 $X\mapsto$ (exp $X)\frak{e}$.  The map $i$ is clearly
 $T$-equivariant open immersion.  Take the Zariski
 closure $\overline{Z^o_S}$ of $i(Z^o_S )$ in $G/B$. Now applying
 \cite{C, Lemma of \S2 } to the
 $T$-stable projective variety $\overline{Z^o_S}\subset
 G/B$, we get that dim$\,Z^o_S=0$ (since $Z^o_S$  does
 not contain any 1-dimensional $T$-stable closed
 irreducible subvarieties). Since any irreducible component of $Z_S$ is
$T$-stable (and closed) in $\frak u^-$ and any closed $T$-stable subset of
$\frak u^-$ contains 0, we get that any irreducible component of $Z_S$ passes
through 0. In particular, $Z_S = \{0\}.$

 Since the variety $X_w$ is Cohen-Macaulay (cf. \cite{Ra
 }, \cite{Ku, Theorem 2.23}, \cite{Ma}), the variety $Y$ is Cohen-Macaulay.
 Assume now that $\# S(w^{-1},v^{-1})= \# S(w,v)= \ell
 (w) = \dim\,Y$, and enumerate the elements of $S(w,v)$
 as $\{ \gamma_1,\dots ,\gamma_\ell  \}$, where $\ell =\ell (w)$.
 By \cite{F, Lemma(a), \S2.4} (since $\dim Z_S=0$),
 the elements $\{ f_{\gamma_j}^Y
 \}_{1\leq j\leq \ell }$ considered as elements of the
 local ring $\Cal O_{0,Y}$ form a regular sequence in $\Cal
 O_{0,Y}$.  Let $I$ be the ideal generated by $\{f^Y_{\gamma_j} \}_{1\leq j
\leq \ell}$ inside the local ring $\Cal O_{0,Y}$. Then there exists an integer
$d >0$ such that $\frak m^d \subset I \subset \frak m$ , where $\frak m \subset
\Cal O_{0,Y}$ is the maximal ideal  (since $Z_S=\{0\}$).  Moreover, by \cite{F,
Lemma (b), \S2.4}, the
 canonical ring homomorphism
 $$
 \frac{\Cal O_{0,Y}}{I}
 [X_1,X_2,\dots ,X_\ell ] \overset
 \sim\to{\longrightarrow}
 \sum_{m\geq 0} I^m / I^{m+1},
 \tag1
 $$
 which takes $X_j$ to the image of $f^Y_{\gamma_j}$ in
 $I/I^2$, is an isomorphism.  In particular,
 $$\align
 \text{ch}\,(\text{gr}\,(\Cal O_{\frak{e},v^{-1}X_w})) &=
 \text{ch}\, (\text{gr}\,(\Cal O_{0,Y})) = \text{ch} (\Bbb C[Y])\quad
\text{($Y$ being
 affine)}\\
 &= \text{ch}\, \biggl( \Cal O_{0,Y} / I \biggr)
 \prod_{j=1}^\ell (1-e^{\gamma_j})^{-1}~,~~~ \text{ by ~(1).} \tag2
 \endalign$$
 But since $\Cal O_{0,Y} / I$ corresponds to the
 0-dimensional variety, it is finite dimensional vector space over $\Bbb C$ and
 hence
 $$
 [ \text{ch}\,(\Cal O_{0,Y} /I )] = \dim \left(\Cal
 O_{0,Y} /I\right).
 \tag3
 $$
 By (2) and (3) we get
 $$
 [\text{ch}\,(\text{gr}\,(\Cal O_{\frak{e},v^{-1}X_w}))] =
 (-1)^\ell d\prod^\ell _{j=1} \gamma^{-1}_j,
 \tag4$$
 where $d:=\dim (\Cal O_{0,Y} /I )$. Thus
 $$\align
 [\text{ch}\,(\text{gr}\,(\Cal O_{\frak{v},X_w}))] &=
 (-1)^\ell \ d\prod^\ell _{j=1}( v\gamma _j)^{-1}\\
 &= (-1)^\ell \, d(-1)^{\#\{ \gamma _j: v\gamma _j\in
 \Delta _-\}}
 \prod_{\beta \in S(w^{-1},v^{-1})} \beta ^{-1}\\
 &= (-1)^{\ell (w)-\ell (v)} d\ \prod_{\beta \in
 S(w^{-1},v^{-1})}
 \beta ^{-1}.
  \tag"(5)" \endalign$$
 This proves the proposition.\qed \enddemo

 \vskip1ex
 \flushpar {\bf (5.3)} {\it Remark.}  When the equivalent
 condition as in the above Proposition (5.2) is
 satisfied, $d$ in fact is an integer $>0$ (as is clear
 from the above proof).
\vskip2ex
  We recall the definition of a rationally smooth point
  in a variety $Y$ (cf.\cite{KL, Appendix}).  \vskip1ex

 \flushpar {\bf (5.4) Definition.} A variety $Y$ of $\dim d$ is
 said to be {\sl rationally smooth} if for all $y\in Y$,
 the singular cohomology $H^i(Y,Y\backslash y,\Bbb Q )=0$
 if $i\neq 2d$ and $H^{2d}(Y,Y\backslash y,\Bbb Q )$ is
 one-dimensional.  A {\sl point} $y_0\in Y$ {\sl is said
 to be rationally smooth} if there exists an open (in the
 Zariski topology) rationally smooth neighborhood of
 $y_0\in Y$.

 A smooth point $y_0\in Y$ is clearly rationally smooth.
 \vskip1ex

 The (b)-part of the following  theorem is the main
 result of this paper.
 \proclaim{(5.5) Theorem}  Fix $v\leq w\in W$.
 \flushpar (a) The point $\frak{v}\in X_w$ is rationally
 smooth $\Leftrightarrow$

 For all $v\leq \theta\leq w$, we have
 $$
 c_{w^{-1},\theta^{-1}}= d_\theta (-1) ^{\ell (w)-\ell
 (\theta)}\,
 \ \prod_{\beta \in S(w^{-1},\theta^{-1})}\ \beta ^{-1},
 \tag1
 $$
 for some constants $d_\theta\in \Bbb C$.
 \flushpar (b)  The point $\frak{v}\in X_w$ is smooth
 $\Leftrightarrow$
 $$
 c_{w^{-1},v^{-1}}
 = (-1) ^{\ell (w)-\ell (v)}\,
 \prod_{\beta \in S(w^{-1},v^{-1})}\ \beta ^{-1} . \tag"(2)"
 $$
 \endproclaim

 \demo{Proof} (a)  By \cite{C, Theorem E}, $\frak{v} \in X_w$ is rationally
smooth if and only if for all $ v\leq
 \theta\leq w$, $\# S(w^{-1},\theta^{-1}) =\ell
 (w)$.  By Proposition (5.2), this is equivalent to the
 requirement that for all $v\leq \theta\leq w$,
 $$
 \text{ch}\,(\text{gr}\,(\Cal O_{\theta,X_w}))
 = d_\theta(-1)^{\ell (w)-\ell (\theta)}
 \prod_{\beta \in S(w^{-1},\theta^{-1})} \ \beta ^{-1},
 $$
 for some $d_\theta\in \Bbb C$.  Now the (a)-part follows
 from  Corollary 3.2(b).
\vskip1ex
 \flushpar (b)  The point $\frak{v}\in X_w$ is smooth if
 and only if  the graded algebra gr$(\Cal O_
 {\frak{v},X_w})$ is isomorphic with the symmetric algebra
 $S[\text{gr}_1(\Cal O_{ \frak{v},X_w})]$.  We first
 prove the `$\Rightarrow$' implication:  So assume that
 $\frak{v}\in X_w$ is smooth.  Then
 $$
 \text{ch}\,(\text{gr}\,(\Cal O_{\frak{v},X_w}))
 = \text{ch}\, (S[\text{gr}_1(\Cal O_{\frak{v},X_w})])
 = \prod_{\gamma \in S} (1-e^\gamma )^{-1},
 $$
 if $\text{ch}\,( \text{gr}_1(\Cal O_{\frak{v},X_w})) =
 \sum_{\gamma \in S}e^\gamma $.  It is easy to see that
 $S\subset v\Delta _+ $ and moreover all the weight
 spaces of $\text{gr}_1(\Cal O_{\frak{v},X_w})$ are
 one-dimensional.  In particular,
 $$
 c_{w^{-1},v^{-1}}
 = [\text{ch}\,( \text{gr}\,(\Cal O_{\frak{v},X_w})) ]
 = \prod_{\gamma \in S\subset v\Delta _+} (-\gamma
 )^{-1}.
 \tag3
 $$
 But since $\frak{v}\in X_w$ is smooth, in particular, it
 is rationally smooth.  So by the (a)-part of the
 theorem,
 $$
 c_{w^{-1},v^{-1}}
 = (-1) ^{\ell (w)-\ell (v)} d_v
 \prod_{\beta \in S(w^{-1},v^{-1})} \beta^{-1},
 \tag4
 $$
 for some positive integer $d_v$ (see Remark 5.6(2)).

 Equating (3) and (4), we get
 $$
 d_v \prod_{\gamma \in S}\gamma  = \pm \prod_{\beta \in
 S(w^{-1},v^{-1})} \beta.
 \tag5
 $$

 Let $Q\subset \frak{h}^\ast $ be the root lattice and
 let  $Q_p:= \Bbb F_p\underset{\Bbb Z}\to{\otimes }Q$ be
 the reduction mod $p$ (for any prime $p$) of $Q$, where
 $\Bbb F_p$ is the prime field of order $p$.
 Reducing the equation (5) mod $p$ (for any prime divisor
 $p$ of $d_v$) and observing that no root mod $p$ is 0 in
 $Q_p$, we get that $d_v=1$.  This proves the
 implication `$\Rightarrow$' of the (b)-part.

 Conversely, assume that
 $c_{w^{-1},v^{-1}}
 = (-1) ^{\ell (w)-\ell (v)}
 \prod_{\beta \in S(w^{-1},v^{-1})} \beta ^{-1} .$
 By Corollary 3.2 (b), this gives
 $$
 [\text{ch}\,( \text{gr}\,(\Cal O_{\frak{v},X_w})) ]
 = (-1) ^{\ell (w)-\ell (v)}
 \prod_{\beta \in S(w^{-1},v^{-1})} \beta ^{-1}.
 \tag6
 $$
 By (5) of the proof of Proposition (5.2), we get that
 $$
 [\text{ch}\,( \text{gr}\,(\Cal O_{\frak{v},X_w})) ]
 = (-1) ^{\ell (w)-\ell (v)} \ d\
 \prod_{\beta \in S(w^{-1},v^{-1})} \beta ^{-1},
 \tag7
 $$
 where
 $d=\dim \left( \Cal O_{0,Y} / I \right)$ (the
 notation is as in the proof of Proposition 5.2).

 By comparing (6) and (7), we get that $d=1$, i.e.,  $I$ is
 the maximal ideal of $\Cal O_{0,Y}$.  In particular, by
 (1) of the proof of Proposition (5.2), gr$(\Cal O_{0,Y})$ is
 graded isomorphic with the polynomial ring $\Bbb
 C[X_1,\dots ,X_\ell ]$.  So we get
 that the point $0\in Y$ is smooth, and hence the point
 $\frak{v}\in X_w$ is smooth.  This proves the theorem
 completely. \qed \enddemo

\vskip1ex
\flushpar
{\bf (5.6)}{\it  Remarks.} (1)~ The (a) part of the above theorem can also be
proved immediately  by combining a result of Dyer \cite{Dy, Proposition  \S3
} with a result of Carrell-Peterson [C,
 Theorem
 E], i.e., we can avoid the use of Corollary 3.2(b).  But our proof has
the advantage that a similar argument (as seen above) gives our criterion
for smoothness as in the (b)-part of the above theorem.
 \vskip1ex
 \flushpar
 (2)~ In the case   (a) as  above (i.e. if ${\frak v}
\in X_w$ is rationally smooth), the constants $d_{\theta}$
 are  in fact  positive
integers for any $v\leq \theta \leq w$  (cf. Remark 5.3).
\vskip1ex
\flushpar
(3) ~ There are some examples
of ${\frak v} \in X_w$ (where $X_w$ is even a codimension one Schubert
variety
in $G/B$) such that $c_{w^{-1},v^{-1}}$ satisfies
condition (1) of the above theorem, but ${\frak v}$ is not a rationally
smooth point of $X_w$ (cf.  Remark 7.9(a)).  In particular, to check the
rational
smoothness of a point ${\frak v} \in X_w$, it is  not sufficient (in
general) to
check the validity of condition (1)
only for $\theta=v$.
\vskip1ex
\flushpar
(4) ~ It is a result of V. V. Deodhar  \cite{D} that any rationally smooth
Schubert variety is in fact smooth for $G = \text{SL}(n)$. This result
has
recently been  extended for any simply-laced $G$ by D. Peterson. As is
well known, this result is false in general for non simply-laced $G$.

\vskip6ex
\centerline{{\bf 6. Singular locus of Schubert varieties in
 rank-$\bold
2$ groups }}
\vskip1ex

As an immediate corollary of Theorem (5.5), we obtain
the
following result determining the singular locus of all the
Schubert varieties in the case of any rank two  group.  I believe
it should be well known, but I did not find it explicitly written
down in the literature. We follow the indexing convention as in Bourbaki
[B].
\vskip2ex
\flushpar
{\bf (6.1) Proposition.}  {\it The following is a complete description
of the
singular locus of the Schubert varieties in the case of rank two groups:}

\flushpar
{\bf Case I}.  $G$ of type $ A_2$ : {\it In this case all the six Schubert
varieties are smooth.}

\flushpar
{\bf Case II.} $G$ of type   $C_2$ :  {\it There are, in all, eight
Schubert varieties.  Out of these only} $X_{r_1 r_2 r_1}$ {\it
is singular and it has singular locus} = $X_{r_1}$.

\flushpar
{\bf Case III.} $G$ of type  $G_2$ :  {\it There
are, in all, twelve Schubert varieties}. {\it Following is the
complete list of singular ones and their singular loci :}

\newpage
\centerline{{\it Singular locus}}\vskip2mm
$$
\matrix \format\l&&\quad\l\\
&(1) ~~~&X_{r_1 r_2r_1}  &-& X_{r_1} \\ \vspace{1\jot}
&(2) ~~~~&X_{r_1 r_2r_1 r_2} &-& X_{r_1 r_2} \\ \vspace{1\jot}
&(3)~~~~ &X_{r_2 r_1r_2 r_1} &-& X_{r_2 r_1} \\ \vspace{1\jot}
&(4)~~~~ &X_{r_1 r_2r_1 r_2r_1} &-& X_{r_1 r_2 r_1} \\ \vspace{1\jot}
&(5) ~~~~ &X_{r_2 r_1r_2 r_1 r_2} &-& X_{r_2}
\endmatrix
$$

\demo{Proof}  As is well known, for any rank-2 group $G$, any $\frak{v}\in X_w$
is
rationally smooth.  (This can also be obtained from Theorem 5.5(a) and the
following Lemma 6.2.) In particular, $c_{w^{-1},v^{-1}}$ satisfies
 identity (1) of Theorem (5.5).  Now the proposition follows
immediately by combining Theorem (5.5)(b) and the following
lemma.\qed\enddemo

 The following lemma can be easily proved by a
 straightforward calculation using the definition of the
 elements $x_{r_i}$ in the nil Hecke ring $Q_W$ (cf.
 Definition 3.1(b)).

 \proclaim{(6.2)  Lemma} For any  group $G$ and any simple reflections $r_1,
r_2 \in W$, we
 have the following (as elements of $Q_W$):
 \endproclaim

 \flushpar (a) $\dsize \quad x_{r_1} x_{r_2}  = \frac{1}{\alpha _1}
 \biggl(
 \frac{1}{\alpha _2}
 \delta _e
 -  \frac{1}{\alpha _2}
 \delta _{r_2}
 - \frac{1}{r_1\alpha _2}
 \delta_{r_1}
 + \frac{1}{r_1\alpha _2} \delta_{r_1r_2}  \biggr)$
 $$\multline x_{r_1}x_{r_2}x_{r_1}  = \\
\frac{1}{\alpha _1}
 \biggl(
 \frac{\alpha _2(\alpha _1^\vee)}{\alpha _2(r_1\alpha _2)}
 (\delta_{e} - \delta_{r_1})+ \frac{1}{\alpha _2(r_2\alpha _1)}  (\delta _{r_2}
-
 \delta _{r_2 r_1})
  -  \frac{1}{(r_1\alpha _2)(r_1r_2\alpha _1)}
 (\delta _{r_1r_2} - \delta _{r_1r_2r_1})
   \biggr)\endmultline\tag"(b)"$$
 $$\multline x_{r_1}x_{r_2}x_{r_1}x_{r_2}   = \\
\frac{1}{\alpha _1}
 \biggl(
 \frac{(m-1)}{\alpha _2(r_1\alpha _2)(r_2\alpha _1)}
 (\delta _e-\delta _{r_2}) -
 \frac{(m-1)}{\alpha _2(r_1\alpha _2)(r_1r_2\alpha _1)}
  (\delta _{r_1} - \delta _{r_1 r_2}) +\\
 \frac{1}{\alpha _2 (r_2\alpha _1)(r_2r_1\alpha _2)}
 (\delta _{r_2r_1} - \delta _{r_2r_1r_2}) -
 \frac{1}{(r_1\alpha _2)(r_1r_2\alpha _1)(r_1r_2r_1\alpha
 _2)}
 (\delta _{r_1r_2r_1} - \delta _{r_1r_2r_1r_2}) \biggr)\endmultline\tag"(c)"$$
$$\multline x_{r_1}x_{r_2}x_{r_1}x_{r_2}x_{r_1}  =\\
 \frac{1}{\alpha _1}
 \biggl(
 \frac{(m-1)(2-m)}{\alpha _2(r_1\alpha _2)(r_2\alpha
 _1)(r_1r_2\alpha _1)} (\delta _e-\delta _{r_1})
  +
 \frac{(2-m)\alpha _2 (\alpha _1^\vee)}{\alpha _2(r_2\alpha
 _1)(r_1\alpha _2)(r_2r_1\alpha _2)}
 (\delta _{r_2}-
 \delta _{r_2r_1}) +\\
 \frac{(m-2)\alpha _2 (\alpha _1^\vee)}{\alpha _2(r_1\alpha
 _2) (r_1r_2\alpha _1) (r_1r_2r_1\alpha _2)} (\delta
 _{r_1r_2} - \delta _{r_1r_2r_1})
  +
 \frac{1}{\alpha _2(r_2\alpha _1) (r_2r_1\alpha _2)
 (r_2r_1r_2\alpha _1)}  (\delta _{r_2r_1r_2} - \delta
 _{r_2r_1r_2r_1}) \\-
 \frac{1}{(r_1\alpha _2)  (r_1r_2\alpha _1)
 (r_1r_2r_1\alpha _2)  (r_1r_2r_1r_2\alpha _1)}
 (\delta _{r_1r_2r_1r_2}- \delta _{r_1r_2r_1r_2r_1})  \biggr) ~,
 \endmultline\tag"(d)"$$
where $m:= \alpha_1(\alpha_2^\vee) \alpha_2(\alpha_1^\vee)$.
\vskip5ex
\centerline{{\bf 7. Singularity of codimension one Schubert varieties in
$\bold G \text{{\bf /}}\bold B$}}
\vskip1ex

Let $w_o$ be the longest element of the  Weyl group $W$ (of
$G$).  As is well known, the codimension one  Schubert varieties in
$G/B$ are precisely of the form $X_w~ ,$ where $w=w_o r_i$ for a
simple reflection $r_i$.  In particular, the number of
such  Schubert varieties in $G/B$ is equal to $n:=$ rank
$G$.  We denote the Schubert variety $X_{w_o r_i}~ (1 \leq i \leq
n)$ by $X_i$.
Let $\chi_i \in {\frak h}_{\Bbb Z }^{*}$ be the $i^{\text{th}}
(1 \leq i \leq n)$
fundamental weight, defined by $\chi_i
(\alpha_{j}^{\vee}) = \delta_{i,j}$.
\vskip2ex
\flushpar
\proclaim{ (7.1) Proposition}  Fix any $1 \leq i \leq n$.  Then for any $v
\in W$ such that $v \leq w_o r_i$,
$$c_{r_iw_0,v^{-1}} = [ \text{ch} \text{(gr}~ {\Cal O}_{{\frak v}, X_i})] =
(-1)^{\mid \triangle_+
\mid -
\ell (v)} \frac{1}{\underset{\beta\in\triangle_+}\to{\prod}\beta} (
w_o \chi_i - v \chi_i), \tag"(1)"$$
where $\text{[~~]}$ is as in $\S$ 3.1(a).

\endproclaim

 \demo{Proof} Consider the i$^{\text{th}}$ fundamental
 representation $V(\chi_i)$ (with highest weight
 $\chi_i$) and define the function
 $$
 \varphi =\varphi _{i,v} : \frak{u}^- \to \Bbb C \text{
 by }\varphi (X) =
 \langle \exp X\cdot e_{\chi_i},\bar v^{-1}e^\ast
 _{w_0\chi_i} \rangle, \text{ for } X\in \frak{u}^-;
 $$
 where $\bar v$ is a representative of $v$ in $N(T)$, $e_{\chi_i}$ (resp.
$e_{w_0\chi_i}$) is a non-zero
 vector in $V(\chi_i)$ of weight $\chi_i$ (resp.
 $w_0\chi_i$) and $e^\ast _{w_0\chi_i} \in
 V(\chi_i)^\ast $ is defined by  $e^\ast
 _{w_0\chi_i}(e_{w_0\chi_i})=1$ and $e^\ast
 _{w_0\chi_i}(v_\mu )$, for any weight vector $v_\mu
 \in V(\chi_i)$ of weight $\mu \neq w_0\chi_i$.  Let $Y$
 be the closed subvariety of the affine space
 $\frak{u}^-$ defined as $Y = \exp^{-1} (U^-\frak{e}\cap
 v^{-1}X_i)$ (cf. proof of Proposition 5.2).  It is easy
 to see that $Y\subset \frak u^-$ is defined set-theoretically by the
 vanishing of the function $\varphi : \frak u^-\to  \Bbb
 C$ (use Lemma 7.2).  Moreover $\varphi$ is obtained by restricting the section
$\chi(\bar{v}^{-1}e^\ast_{w_0\chi_i}) \in H^0(G/B, \frak L(\chi_i))$ to
$U^-\frak e$ (and using the identification exp: $\frak u^- \to U^-\frak e
\subset G/B$ ), where $\chi$ is the Borel-Weil homomorphism (cf. Proof of Lemma
2.4). But the line bundle  $\frak L(\chi_i)$ on $G/B$ corresponds to the
irreducible divisor $X_i \subset
G/B$ with multiplicity $1$ (use, e.g.,  the Chern class calculation for the
line bundle $\frak L(\chi_i) $ ). This, in particular, implies that the ideal
$I$ of the irreducible hypersurface $Y\subset \frak u^-$  (with the reduced
structure) is generated by the function $\varphi$ (cf. also \cite{C2,
Proposition 4.6}. This gives that  (as graded $T$-algebras),
 $$
 \text{gr}\,(\Cal O_{\frak{e},v^{-1}X_i})
 \approx
S(\frak{u}^{-^\ast}) /\langle [\varphi ] \rangle
 ,
 \tag2
 $$
 where (as earlier) $S(\frak{u}^{-^\ast })$ is the symmetric algebra
 of $\frak{u}^{-^\ast }$ and $\langle [\varphi ] \rangle$
 denotes the (homogeneous) ideal generated by the least
 degree non-zero homogeneous component $[\varphi ]$ of
 $\varphi$.  From the definition of $\varphi$, it is easy
 to see that $[\varphi ]$ is a weight vector for the
 adjoint action of $T$ on $\frak{u}^-$  with  weight
 $\chi_i - v^{-1}w_0\chi_i$.  So by (2),
 $$\gather
 \text{ch}\, (\text{gr}\, \Cal O_{\frak{e},v^{-1}X_i})
 =
 (1-e^{\chi_i- v^{-1}w_0\chi_i} )
 \prod_{\beta \in\Delta _+} \ (1-e^\beta)^{-1},\\
 \intertext{and hence}
 [\text{ch}\, (\text{gr}\, \Cal O_{\frak{e},v^{-1}X_i}) ]
 =
 (-1) ^{\#\Delta _+ }~
 \frac{(v^{-1}w_0\chi_i-\chi_i)}{\underset{\beta \in \Delta
 _+}\to{\prod}\beta }.
 \tag3\endgather$$
 (Observe that by Lemma (7.2), $v^{-1}w_0\chi_i-\chi_i\neq 0$, since by
 assumption $v\leq w_0r_i$.)  By applying $v$ to (3) we
 get
 $$
 [\text{ch}\, (\text{gr}\, \Cal O_{\frak{v},X_i}) ]
 =  (-1)^{\# \Delta_+ -\ell (v)}
 \frac{(w_0\chi_i-v\chi_i)}{\underset{\beta \in \Delta
 _+}\to{\prod} \beta }.
 $$
 This proves the second equality of (1).  First equality of (1) of course
follows from Corollary 3.2(b). \qed

\enddemo

 \proclaim{(7.2) Lemma} For any simple reflection $r_i$
 and any $v\in W$, $v\leq w_o r _i$ if and only if
 $\chi_i\neq v^{-1}w_0\chi_i$.
 \endproclaim

 \demo{Proof}  Let $Z\subset G$ be the zero set of the
 function $\hat{\varphi} : G\to \Bbb C$ given by
 $\hat{\varphi} (g)= \langle ge_{\chi_i},e^\ast
 _{w_0\chi_i} \rangle$ (where $e_{\chi_i}$  and $e^\ast
 _{w_0\chi_i}$ are as in the proof of Proposition 7.1).  Then
 clearly $Z$ is $B$-stable under the left as well as
 right multiplication.  In particular, $Z/B= \cup X_j$,
 where $j$ runs over some subset $S\subset \{ 1,\dots
 ,n\}$.  Clearly $i\in S$, whereas for $j\neq i,\ j\notin
 S$, and hence $Z/B=X_i$.  Hence $v\leq w_0r
 _i\Leftrightarrow v\in X_i=Z \Leftrightarrow v\chi_i
 \neq w_0\chi_i$.\qed

 \enddemo

 \proclaim{(7.3) Lemma} Assume that $v\leq w_0r _i$.
 Then $\chi_i -v^{-1}w_0\chi_i$ is multiple of a root
 $\beta $ if and only if $\pm v\beta \notin S(r _i
 w_0,v^{-1})$.  In particular, $\chi_i -v^{-1}w_0\chi_i$
 is multiple of a root if and only if  \#$S(r
 _iw_0,v^{-1})= N-1$, where $N := \# \Delta _+$.
 \endproclaim

 \demo{Proof}  If $\pm v\beta \notin S(r
 _iw_0,v^{-1})$, then by the above Lemma (7.2), $r
 _\beta v^{-1}w_0 \chi_i=\chi_i$.  In particular,
 $\chi_i-v^{-1}w_0\chi_i$ is a multiple of $\beta $.

 Conversely, assume that
 $$
 \chi_i-v^{-1}w_0\chi_i=n\beta,
 \tag1$$
 for some number $n$ and $\beta \in \Delta $ .  By
 Lemma (7.2), $n\neq 0$.  To prove that $\pm v\beta
 \notin S(r _iw_0,v^{-1})$, it suffices to show
 (again by Lemma 7.2) that $r _\beta v^{-1}w_0
 \chi_i = \chi_i$:  By (1),
 $$\gather
 \langle \chi_i - v^{-1} w_0 \chi_i,\beta^\vee \rangle =
 2n, \quad \text{and}\tag2\\
 \langle \chi_i + v^{-1} w_0 \chi_i,\beta^\vee \rangle =
 \frac{2}{n\langle \beta ,\beta  \rangle}
 \langle \chi_i + v^{-1} w_0 \chi_i , \chi_i - v^{-1} w_0
 \chi_i \rangle
=0.
 \tag3\endgather
 $$
 Combining (2) and (3) we get $\langle
 -v^{-1}w_0\chi_i,\beta ^\vee \rangle = n$; and hence
 $r _\beta v^{-1}w_0 \chi_i := v^{-1}w_0\chi_i -
 \langle v^{-1}w_0\chi_i,\beta ^\vee \rangle\beta =
 v^{-1}w_0\chi_i+n\beta =\chi_i$ (by (1)).

 The `in particular' statement of the lemma follows from
 Deodhar's conjecture (cf. Theorem 5.1).\qed
 \enddemo

By virtue of  Proposition (7.1),  Lemma (7.3),  and Theorem 5.5(b), we get the
following
characterization of the smooth points in  the Schubert varieties $X_i$.

\proclaim{(7.4) Proposition } Let $X_i ~~ (1 \leq i \leq n)$ be a
codimension one  Schubert variety.  Then, for any $v \leq
w_o r_i \in W$, the following are equivalent:
\roster \item"(a$_1$)"
{}~ ${\frak v} \in X_i$ is smooth.
\item"(a$_2$)"  ~ $c_{r_i w_o,v^{-1}} = (-1)^{N-1-\ell (v)}
\frac{1}{\beta_1 \cdots \beta_{N-1}}$, for some positive roots
$\{\beta_1, \cdots, \beta_{N-1} \}$
(where $N=$ dim $G/B$).
\item"(a$_3$)"~ $\chi_i- v^{-1} w_o \chi_i$ is a root. \endroster

   $\ \ \ $In particular, $X_i$ is smooth if and only if
  $\chi_i-w_0\chi_i$ is a root.
  \endproclaim
\flushpar
{\bf (7.5)}
{\it  Remark.} If ${\frak v} \in X_i$ is smooth, then the set
$\{\beta_1, \cdots, \beta_{N-1}\}$, as in ($a_2$) above, coincides
with the set $S(r_{i}w_{0},v^{-1})$ (by Theorem 5.5 (b)).

  \demo{Proof (of Proposition 7.4)} As follows from
  Theorem 5.5(b), (a$_1$)$\Rightarrow$(a$_2)$.  The
  implication (a$_2$)$\Rightarrow$(a$_3)$ follows from
Proposition
  (7.1).  So we come to the proof of (a$_3)\Rightarrow$
  (a$_1)$:

  By Theorem 5.5(b), we need to show that
  $$
  c_{r _iw_0,v^{-1}} = (-1)^{N-1-\ell (v)}
  \prod_{\beta \in S(r _iw_0,v^{-1})} \beta^{-1} .
  \tag1$$

  By (a$_3$), $\gamma := v\chi_i-w_0\chi_i$ is a root
 (and
  in fact is positive since $v\leq w_0$).  In particular,
  by Proposition (7.1),
  $$
  c_{r _iw_0,v^{-1}} =
  (-1)^{N-1-\ell (v)}
  \frac{\gamma}{\prod_{\beta \in \Delta _+} \beta }~.
  \tag2$$
  But by Lemma (7.3), $S(r _iw_0,v^{-1}) = \Delta
  _+\backslash \{ \gamma \}$, and hence (1) follows from
  (2).  This proves the implication (a$_3)\Rightarrow
  $(a$_1)$.

  The `in particular' statement of the proposition
 follows
  from the equivalence of (a$_1$) and (a$_3$) since
 $X_i$
  is smooth if and only if $\frak{e}\in X_i$ is
  smooth.\qed
  \enddemo
  By the same proof as above for the implication
  (a$_3$)$\Rightarrow$(a$_1)$ (alternatively, by using
 Lemma
  (7.3) with \cite{C, Theorem E}) we obtain the
 following:

  \proclaim{(7.6) Corollary} With the notation as in
  Proposition (7.4), $\frak{v}\in X_i$ is rationally
  smooth if and only if for all $v\leq \theta \leq
  w_0r _i$, $\chi_i-\theta^{-1} w_0\chi_i$ is
  multiple of a root $\beta _\theta$ (depending upon
  $\theta$). \qed
  \endproclaim

  We follow the indexing convention of simple roots as in
  \cite{B, Planche I-IX}.  The following  lemma follows
  easily from the explicit knowledge of roots, coroots,
  fundamental weights etc. as given in loc. cit.

  \proclaim{(7.7) Lemma}Let $G$ be a simple algebraic
  group.  Then for any fundamental weight $\chi_i (1\leq
  i\leq n)$,
  \flushpar (a) $\ \ \ \chi_i-w_0 \chi_i$ is a (positive)
 root
  precisely in the following cases ($A_n$ etc. denotes
 \par the
  type of $G$):
    $$\matrix \format\l&&\quad\l\\  (\text{{\text a}}_1)
 &A_n  &(n\geq 1)  &; &i=1,n \\
     (\text{{\text a}}_2) &C_n  &(n\geq 2) &; &i=1.
 \endmatrix
 $$
  \flushpar (b) $\ \ \ \chi_i-w_0\chi_i$ is multiple of a
 root but
  not a root itself, precisely in the following \par cases:
  $$
 \matrix \format\l&&\quad\l\\
   (\text{{\rm b}}_1 )  &B_n &(n\geq 3)  &;\ &  i=1,2\\
     (\text{{\rm b}}_2 )  &C_n &(n\geq 2)  &;\ & i=2\\
     (\text{{\rm b}}_3 )  &D_n &(n\geq 4)  &;\ & i=2 \\
     (\text{{\rm b}}_4 )  &E_6           & &;\ & i=2 \\
     (\text{{\rm b}}_5 )  &E_7           & &;\ & i=1 \\
     (\text{{\rm b}}_6 )  &E_8           & &;\ & i=8 \\
     (\text{{\rm b}}_7 )  &F_4           & &;\ & i=1,4\\
     (\text{{\rm b}}_8 )  &G_2            & &;\ & i=1,2.
 \endmatrix
 $$
  \endproclaim

\vskip1ex
 As a consequence of the above lemma, we get the
 following complete list of codimension-1 Schubert
 varieties which are smooth or rationally smooth.

 We assume that $G$ is a simple group in the following
 proposition.

 \proclaim{(7.8) Proposition}  {\rm (c)} The following is
 a complete list of codimension one Schubert varieties
 $X_i$ which are smooth:
 $$\matrix\format \l&&\quad\l\\
 (\text{{\rm c}}_1 )  &A_n &(n\geq 1)&: &i=1,n\\
   (\text{{\rm c}}_2 )&C_n &(n\geq 2)&: &i=1.
 \endmatrix
 $$
 \flushpar (d)  The following is a complete list of
 codimension one Schubert varieties $X_i$ which are
 rationally smooth but not smooth:
 \roster
   \item"(d$_1$)"  $C_2 : i=2$
   \item"(d$_2$)"  $G_2 : i=1,2$
 \item"(d$_3$)" $B_n$ ($n\geq 3): i=1$.  \endroster
 \endproclaim

 \demo{Proof} The (c)-part follows immediately by
 combining Proposition (7.4) with Lemma (7.7).

 To prove the (d)-part, in view of Corollary (7.6) and
 Lemma (7.7), it suffices to show that in all the cases
 covered by (b) of Lemma (7.7) but not in the list (d)
 above, there exists a $\theta\in W$ such that
 $\chi_i-\theta^{-1}w_0\chi_i$ is not a multiple of any
 root (by Lemma 7.2, such a $\theta$ will automatically
 satisfy $\theta\leq w_0r i$), whereas in the cases
 covered by (d), $\chi_i-\theta^{-1}w_0\chi_i$ is indeed
 multiple of a root for any $\theta\in W$:

 We freely use the notation without explanation from
 \cite{B; Planche I-IX}.  In the cases $( B_{n\geq 3};\
 i=2 ),\ ( C_{n\geq 3};\ i=2 )$ and $( D_{n\geq 4};\ i=2
 )$ take any $\theta\in W$ satisfying $\theta(\epsilon_1)=\epsilon_1$,\
 $\theta(\epsilon_3)=\epsilon_2$.
 Then $\chi_i-\theta^{-1}w_0\chi_i$ is not a multiple of
 any root.

In the cases ($E_6$; $i=2$), ($E_7$; $i=1$),
 and ($E_8$; $i=8$), $\chi_i$ is the highest root
 $\alpha_0$.  In these cases take any $\theta\in W$
 satisfying $\theta(\alpha_2)= \alpha_0$  (observe that
 $-w_0\alpha_0=\alpha_0$ and the $W$-orbit $W\cdot
 \alpha_0$ consists of all the roots), then
 $\chi_i-\theta^{-1}w_0\chi_i$ is not a multiple of any
 root.

 In the case ($F_4$; $i=1$), (resp. $F_4$; $i=4$),
 $\chi_1$ (resp. $\chi_4$) is the highest (resp. a short)
 root, in particular, $W\cdot \chi_1$ consists of all the
 long (resp. short) roots.  Take any $\theta\in W$ satisfying
 $\theta(\epsilon_2+\epsilon_3)=\chi_1$ (resp. $\theta\left(
 \frac{\epsilon_1+\epsilon_2+\epsilon_3+\epsilon_4}{2}\right) =\chi_4 $), then
 $\chi_i-\theta^{-1}w_0\chi_i$ is not a multiple of any
 root.

 For (C$_2$; $i=2$) and  (G$_2$; $i=1,2$), it is easy to
 see that $\chi_i-\theta^{-1}w_0\chi_i$ is  multiple of
 a root, for all $\theta\in W$.

So finally we come to
 ($B_{n\geq 3};\ i=1$):  In this case, $-w_0=\text{Id.},\
 \chi_1=\epsilon_1$ (a short root), and hence $W\cdot \chi_1=\{
 \pm \epsilon_i\}_{1\leq i\leq n}$.  In particular,
 $\chi_1-\theta^{-1}w_0\chi_1$ is multiple of a root for
 all $\theta\in W$.  This finishes the proof of the
 (d)-part of the proposition. \qed

 \enddemo

 \flushpar {\bf (7.9)}{\it Remarks.}  (a)  In all the cases
 covered by Lemma 7.7(b) but not contained in Proposition
 7.8(d),  identity (1) of Theorem 5.5(a) is satisfied
 for $w=w_0r _i$ and $\theta=e$ but is violated for
 some $e\leq \theta\leq w$ (use Proposition 7.1, Lemma 7.3 and Theorem 5.5(a)).

 \flushpar (b)  I am informed that the (c) part of the
 above proposition, as well as the equivalence  of ($a_1$)
 and ($a_3$) in Proposition (7.4) for $v=e$ was contained in
 an earlier longer version of \cite{C} (cf. \cite{C2, \S4}).  Of course
($d_1$), ($d_2$)
 are  very well known, and  example ($d_3$) was known to
 be rationally smooth due to  Boe \cite{Bo}.
\vskip5ex

 \centerline{{\bf 8.  Extension of results to the
 Kac-Moody case}}
 \vskip2ex
 \flushpar {\bf (8.1) Notation.}  We will follow the
 notation (often without explaining) from \cite{Ku; \S
 1}:
 In particular, throughout this section $G=G(A)$ denotes the
 complex Kac-Moody group associated to an arbitrary
 $n \times n $  generalized Cartan matrix $A$
 (we do not put symmetrizability restriction on $A$),
 with  the standard Borel subgroup $B$, and the
 standard maximal torus $T\subset B$.  There is a
  Weyl group $W\approx N(T)/T$ associated to the pair
 $(G,T)$ (where $N(T)$ is the normalizer of $T$ in $G$).
  The Weyl group
 $W$ is a  Coxeter group with the   simple
 reflections $\{ r_i \}_{1\leq i\leq n }$ as Coxeter
 generators ($r_i$ is nothing but the reflection through
 the simple root $\alpha_i$).  Hence, for any $w\in W$,
 we can talk of its  length $\ell (w)$ and also have
  Bruhat partial ordering $\leq$ in $W$.

 The  Kac-Moody algebra  $\frak{g} =
 \frak{g}(A)$ admits the root space
 decomposition:
 $$
 \frak{g}=\frak{h} \oplus \sum_{\alpha\in
 \Delta _+ \subset \frak h^\ast}
 (\frak{g}_\alpha \oplus \frak{g}_{-\alpha}),
 $$
 where $\frak{g}_\alpha := \{ X\in \frak{g}:
 [h,X]=\alpha(h)X,\ \text{for all }h\in \frak{h} \}
 $ is the $\alpha$-th root space, $\frak{h} :=$  Lie $T$ is the
  standard Cartan subalgebra of $\frak{g}$,
  and $\Delta _+ := \{
 \alpha \neq 0 \in \sum^n _{i=1} \Bbb Z_+ \alpha_i
 : \frak{g}_\alpha \neq 0\}$ is the set
 of   positive roots.  We denote $\Delta _- =
 -\Delta _+$ and $\Delta := \Delta _+\cup \Delta _-$.
 The Weyl group $W$ preserves $\Delta$.  The  set of
 real roots $\Delta ^{\text{re}}\subset \Delta $ is
 defined to be $W$.  $\{ \alpha_1,\dots ,\alpha_n  \}$
 and the  set of imaginary roots $\Delta
 ^{\text{im}} := \Delta \backslash \Delta ^{\text{re}}$.
 We set $\Delta ^{\text{re}}_+ = \Delta _+\cap \Delta
 ^{\text{re}}$ (resp. $\Delta ^{\text{re}}_- = \Delta_- \cap
 \Delta ^{\text{re}}$); $\Delta ^{\text{im}}_+$ and
 $\Delta ^{\text{im}}_-$ have similar meanings.  We denote
 by $\tilde{\Delta}_+$  (resp.
 $\tilde{\Delta}_-$) the indexed set of
 positive (resp. negative) roots such that each root
 occurs exactly as many times as the dimension of the corresponding root
 space. Recall that the real root spaces are of dimension
 one.

 The group $G$ (in particular the torus $T$ ) acts on $G/B$ by the left
multiplication.
For any $w\in W$, the Schubert variety $X_w$ is by
 definition the closure of $B\bar{w}B/B$ in $G/B$, where
$\bar{w}$ is a preimage of $w$ in $N(T)$ and $G/B$  is endowed with the Zariski
topology as in
 \cite{S}.  Of course, $X_w$ is $T$-stable.  As is well
 known (by the Bruhat decomposition), $X_w= \underset
 v\leq w\to{\cup} B\bar{v} B/B$.  In particular, for any $v\leq
 w,\ \frak v:=\bar{v}B\in X_w$ and it is a $T$-fixed point.  We will
 always endow $X_w$ with the  stable variety
 structure as given in \cite{Ku; \S1}.  With this
 structure $X_w$ is an irreducible projective variety of
 dim $\ell (w)$.

 For any real root $\beta $, there exists a unique
 additive one-parameter subgroup $U_\beta $ and a
 homomorphism $u_\beta :\Bbb C \to G$ satisfying $u_\beta (\Bbb C)= U_\beta $
and such that
 $$
 t u_\beta (z) t^{-1} =
 u_\beta (e^\beta (t)z),  $$
 for any $z \in \Bbb
 C$,
and $t\in T$.
 Furthermore, for any $w\in W$,
 $
 \bar{w} U_\beta \bar{w}^{-1} =  U _{w\beta
 }. $

 Now let $U^-$ be the subgroup of $G$ generated by the
 one-parameter groups $\{ U_\beta \}_{\beta \in \Delta
 _-^{\text{re}}}$.  Then the map $U^- \to G/B$, taking
 $g\mapsto g\frak{e}$ is injective and moreover
 $U^-\frak{e}\subset G/B$ is an open subset.

 For any $\lambda\in \frak h^\ast_{\Bbb Z}$, recall the definition of the
 line bundle $\frak{L}(\lambda):= G\underset
 B \to{\times} \Bbb C_{-\lambda}\to G/B$ from
 \cite{Ku; \S2.2}.  For dominant $\lambda\in   \frak h^\ast_{\Bbb Z}$, let
 $V^{\max}(\lambda)$ be the maximal integrable highest
 weight $G$-module with highest weight $\lambda$ (cf.
 \cite{Ku, \S1.5}, where it is denoted by
 $L^{\max}(\lambda)$).  Define
 $$
 H^0(G/B,\frak{L}(\lambda)) =
 \text{Inv limit}_{w\in W}
 H^0(X_w,\frak{L}(\lambda)_{\vert_{X_w}}) .
 $$
 The highest weight space
 $\Bbb C_\lambda :=V^{\text{max}}(\lambda)_{(\lambda)}$ of
 $V^{\text{max}}(\lambda)$ is one dimensional.  Define
 the map
 $$
 \chi = \chi_\lambda: V^{\text{max}}(\lambda)^\ast
 \longrightarrow
 H^0(G/B,\frak{L}(\lambda))
 $$
 by $\chi(f) (gB) = (g,  (g^{-1}f)\vert{\Bbb C}_\lambda) \text{ mod}\,B$,
 for $f\in V^{\text{max}}(\lambda)^\ast$, and $g\in G$.

 The following  result is due to Kumar \cite{Ku, Theorem 2.16}
 (and also  Mathieu \cite{Ma}).

 \proclaim{(8.2) Theorem}  The map $\chi_\lambda$ as above
 is an isomorphism.  Moreover, for any $v\leq w\in W$, it
 induces an isomorphism
 $$
 \chi_\lambda(v,w) :
 (v^{-1}V_w^{\text{max}}(\lambda))^\ast
 \tilde\rightarrow
 H^0(v^{-1}X_w,
 \frak{L}(\lambda)_{\vert_{v^{-1}X_w}}),
 $$
 making the following diagram commutative:
 $$\CD
 V^{\text{max}} (\lambda)^\ast
 @>\chi_\lambda>>
 H^0(G/B,\frak{L}(\lambda)) \\
 @VVV  @VVV\\
 (v^{-1} V^{\text{max}}_w (\lambda))^\ast
 @>\chi_\lambda(v,w)>>
 H^0(v^{-1}X_w,
 \frak{L}(\lambda)_{\vert_{v^{-1}X_w}})
 \endCD
 $$
 where $V^{\text{max}}_w (\lambda)\subset V^{\text{max}}
 (\lambda)$ is the $B$-submodule generated by the extremal
 weight space $V^{\text{max}} (\lambda)_{(w\lambda)}$ of weight $w\lambda$, and
the vertical maps are the canonical restriction maps.
 \endproclaim

 For any non-zero $e_\lambda\in
  \Bbb C_\lambda$, define
 $e_\lambda^\ast \in V^{\text{max}} (\lambda)^\ast$ as
 $e_\lambda^\ast(e_\lambda) =1$ and $e_\lambda^\ast(y)=0$, for any
 weight-vector $y$ of weight $\mu \neq \lambda$.   Now define the section
 $s_{e_\lambda} \in H^0(G/B,\ \frak{L}(\lambda))$
 by $s_{e_\lambda} = \chi_\lambda (e_\lambda^\ast)$.
\vskip2ex

 The following  lemma follows immediately from the
 Birkhoff decomposition \cite{KP, \S3}.
 \proclaim{(8.3) Lemma}  The zero set of $s_{e_\lambda} , Z(s_{e_\lambda})
  = G/B\backslash (U^- \frak{e})$,
 if $\lambda \in D^o$, where (as in \S1) $D^o$ is the set  of dominant regular
weights.\qed
 \endproclaim

 The line bundle
 $\frak{L}(\lambda)_{\vert_{v^{-1}X_w}}$ on the
 projective variety $v^{-1}X_w$ is ample for any $v\leq
 w\in W$ and $\lambda \in D^o$.  In
 particular, by Lemmas (2.3) and (8.3),
 $U^-\frak{e}\cap v^{-1}X_w$ is an affine open
 subset of $v^{-1}X_w$.

 Define the $T$-equivariant map (cf. \S2.6)
 $$
  \varphi _\lambda(v,w) :
 (v^{-1} V^{\text{max}}_w (\lambda))^\ast
 \otimes
  \Bbb C_\lambda \rightarrow
 \Bbb C [  U^- \frak{e}  \cap v^{-1}X_w]
 ~\text{by} $$
 $$(\varphi _\lambda (v,w) (f\otimes e_\lambda)) (x) s_{e_\lambda }
 (x)
 =
 (\chi_\lambda(v,w) f) (x) ,
  $$
 for $f\in (v^{-1} V^{\text{max}}_w (\lambda))^\ast ,\
 e_\lambda \neq 0\in \Bbb C_\lambda$ and $x\in U^-\frak{e}\cap v^{-1}
 X_w $.  (We set $\varphi _\lambda (v,w) (f \otimes 0) =0.$) By Lemma (8.3),
the map $\varphi _\lambda (v,w)$
 is well defined, and is injective by Theorem (8.2).
 Moreover, as in \S2.7, for any $\lambda\in D^o$ and $\mu
 \in D$, the following diagram is commutative:
 $$
   \gather
   (v^{-1} V^{\text{max}}_w (\lambda))^\ast \otimes
 \Bbb C_\lambda
  \quad
   \overset{\delta_{\lambda,\mu }(v,w)}\to\hookrightarrow
  \quad
   (v^{-1} V^{\text{max}}_w (\lambda+\mu ))^\ast \otimes
 \Bbb C_{\lambda+\mu} \\ \vspace{4\jot}
   \sideset{\varphi_\lambda (v,w)}\and\to\searrow
  \qquad\qquad\qquad
  \sideset\and{\varphi_{\lambda +\mu }(v,w)}\to\swarrow\\
 \vspace{2\jot}
   \Bbb C [U^- \frak{e} \cap v^{-1} X_w],
   \endgather
  $$
 where the map $\delta_{\lambda,\mu}(v,w)$ is defined
 as in Lemma (2.7).  Taking the limit of the
 maps $\varphi _\lambda(v,w)$, we get the $T$-equivariant
 map
 $$
 \varphi (v,w) : \underset \lambda\in D^o\to{\text{limit}}
 \rightarrow
 ((v^{-1}V^{\text{max}}_w(\lambda))^\ast \otimes \Bbb
 C_\lambda)
 \rightarrow \Bbb C [U^-\frak{e} \cap v^{-1}X_w].
 $$
 The following proposition follows easily from Lemma
 (2.3) and Theorem (8.2).

 \proclaim{(8.4) Proposition}  The above map $\varphi
 (v,w)$ is an isomorphism for any $v\leq w\in W$.\qed
 \endproclaim

 Define the Lie subalgebra  $\frak{u}^- =
 \underset{\alpha\in \Delta _-}\to{\oplus}
 \frak{g}_{\alpha} $ of $\frak{g}$ and (for any
 $m>0$) the ideal $\frak{u}^- _m$ of
 $\frak{u}^-$ by
 $$
 \frak{u}^- _m =
 \oplus \Sb \alpha\in \Delta _-\\\vert\alpha \vert\geq m \endSb
 ~\frak{g}_{\alpha},
 $$
 where for a root
 $$
 \alpha= \sum m_i\alpha_i,\  \vert \alpha \vert:= \vert \sum m_i
 \vert .
 $$

 The quotient algebra $F_m(\frak{u}^-) :=
 \frak{u}^-/\frak{u}^-_m$ is a finite dimensional nilpotent
 algebra.  Let $F_m(U^-)$ be the associated unipotent
 complex algebraic group.  Corresponding to the Lie
 algebra homomorphism $\frak{u}^-\to
 F_m(\frak{u}^-)$, there is associated a group
 homomorphism $\theta_m: U^-\to F_m(U^-)$. We state the following simple lemma
without proof.
 \proclaim{(8.5) Lemma} Fix $v\leq w\in W$.  Then there
 exists a positive number $m_0(v,w)$ such that
 $$
 \theta_m(v,w) : U^-\frak{e} \cap v^{-1}X_w \to
 F_m(U^-)
 $$
 (got by restricting the map $\theta_m$) is a closed
 immersion for all $m\geq m_0(v,w)$.
 \endproclaim

 By an argument identical to the proof of Theorem (2.2)
 (as given in \S2.12), and Corollaries (3.2) (using
 Proposition  8.4, Lemma  8.5,  and \cite{Ku, Theorem 3.4}) we get the
following
 analog of Theorem (2.2) and Corollaries (3.2) for an
 arbitrary Kac-Moody group $G$.

 \proclaim{(8.6) Theorem}  Let $G$ be an arbitrary Kac-Moody group.
\roster
 \item "{\rm (a)}"   For  any $v \leq w \in W$,  $\text{gr }
{\Cal O}_{{\frak v},X_w}$  is an admissible $T-$module and  moreover
$$\text{ch~(~gr}~ {\Cal O}_{{\frak v},X_w})= * b_{w^{-1},v^{-1}},$$
as elements of $\widetilde{Q(T)}$ .
 \item "{\rm (b)}" For any $v\leq w  \in W$, $b_{w^{-1},v^{-1}}\neq 0$ if and
only if $v\leq
 w$, and in this case it has a pole of order exactly
 equal to $\ell (w)$. Further, there exist $\beta_1, \dots, \beta_N \in
\tilde\Delta^+$ (for some $N >0$) such that
 $$
 \biggl( \prod_{j=1}^N (1-e^{\beta _j}) \biggr)
 b_{w^{-1},v^{-1}}
 \in R(T).
 $$
 \item "{\rm (c)}" $[\ast b_{w^{-1},v^{-1}} ] =
 c_{w^{-1},v^{-1}}$; and hence for any $v\leq w,\
 [\text{{\rm ch}}\,(\text{{\rm gr}}\,\Cal O_{\frak v,X_w})] =
 c_{w^{-1},v^{-1}} $, as elements of $Q(\frak{h})$.

{\it In particular, $c_{w, v} \neq 0$ if and only if} $v
\leq w$. \qed
\endroster
 \endproclaim

 We extend Proposition (5.2) to the Kac-Moody case.

 \proclaim{(8.7) Proposition}  Let  $G$ be an arbitrary Kac-Moody group and let
$v\leq w\in W$.  Then
 $$
 \sharp S_{w^{-1},v^{-1}}  \,= \ell (w)
 \Leftrightarrow [ \text{ch} \, ( \text{gr} \, \Cal
 O_{\frak{v},X_w})]
 = d(-1)^{\ell (w)-\ell (v)}
 \prod_{\beta\in S(w^{-1},v^{-1})} \ \beta ^{-1},
 $$
 for some $d\in \Bbb C$; where $S(w^{-1},v^{-1})= \{ \alpha\in \Delta
_+^{\text{re}}: v^{-1}r
 _\alpha\leq w^{-1}\}$.
 \endproclaim

 \demo{Proof} The proof is very similar to the proof of
 Proposition (5.2).  But we need to make the following
 modifications:

 Define $Y= U^- \frak{e}\cap v^{-1}X_w$.  Fix any
 regular $\lambda\in D^o$ and a highest weight vector
 $e_\lambda\in V^{\text{max}}(\lambda)$ and consider the
 element $e_\lambda^\ast\in V^{\text{max}}(\lambda)^\ast $ as
 in \S8.2.  For any root $\alpha\in \Delta
 _+^{\text{re}}$, choose a non-zero root vector
 $X_{\alpha}\in \frak{g}_\alpha$ and define the map
 $\theta_\alpha: U^- \to \Bbb C$ by
 $
 \theta_\alpha(g) = e_\lambda^\ast(X_\alpha ge_\lambda), ~
 \text{for }g\in U^-.
 $
 We claim that $\theta_\alpha(g) \neq 0$, for any $g\neq e\in
 U_{-\alpha}$:

 Write $g=\exp(zX_{-\alpha})$, for some $z\neq 0 \in \Bbb C$ ;
 where $X_{-\alpha}$ is the root vector corresponding to
 the (real) root $-\alpha$ such that
 $[X_\alpha,X_{-\alpha}] = \alpha^\vee$ (cf. \cite{K, exercise 5.1}).
 Then $$\align
  \theta_\alpha (g) &=
  e^\ast_\lambda(X_\alpha \exp(zX_{-\alpha})e_\lambda) \\
 &=  e^\ast_\lambda (zX_\alpha X_{-\alpha} e_\lambda) \\
 &=  e^\ast_\lambda (z [X_\alpha,X_{-\alpha}] e_\lambda)\\
 &= z \langle \lambda,\alpha^\vee \rangle\\
 &\neq 0 ~,\qquad \text{since } \lambda \text{ is regular.}
 \endalign
 $$
 Idetifying $U^-
\simeq U^-\frak e$, we can (and do) consider $\theta_\alpha$ as a function on
$Y$.
 Now define
 $$
 Z_S =
 \{ x\in Y :
  \theta_\alpha(x) =0, \text{ for all } \alpha\in S:= S(w,v)\}.
 $$

 Rest of the argument to prove the proposition is similar
 to the proof of Proposition (5.2) provided we replace
 $\frak{u}^-$ by $U^-$ and use the following  simple
\enddemo

 \proclaim{(8.8) Lemma}  For any $v\leq w\in W$, one
 dimensional $T$-orbits in  $U^-\frak{e}\cap
 v^{-1}X_w$ are precisely of the form $(U_{-\beta
 }\backslash e)\frak{e}$, where $\beta $ ranges over
  (positive real) roots $\in S(w,v)$.
 \endproclaim

 \demo{Proof}  By the Bruhat decomposition
 $$
 X_w = \underset {\theta\leq w}\to{\cup} U
 \theta\frak{e} =
 \underset {\theta\leq w}\to{\cup}
 \theta(\theta^{-1}U\theta\cap U^-) \frak{e},
 $$
 one-dimensional $T$-orbits contained in
 $v^{-1}X_w$ are precisely of the form $I_{\theta,\beta
 }:= v^{-1}\theta (U_{-\beta }\setminus e)\frak{e}$, where
 $\theta\leq w$ and $\beta \in \Delta _+\cap
 \theta^{-1}\Delta _-$.  (We are using the fact that any root in  $ \Delta
_+\cap
 \theta^{-1}\Delta _-$ is a real root and moreover for any real root $\beta$,
$d\beta$ is not a root for any $d>1$.) If $v=\theta$, clearly
 $I_{\theta,\beta }\subset U^-\frak{e}$, and
 moreover $\beta \in \Delta _+\cap v^{-1}\Delta _-
 \Leftrightarrow \beta \in \Delta _+^{\text{re}}$ and
 $vr_\beta < v$ (by \cite{BGG, Corollary 2.3}).  So assume
 that $v\neq \theta$.  By Bruhat decomposition for $SL(2)$, we get
$$(U_{-\beta} \frak e) \cup \{r_\beta \frak e\} = \overline{Br_\beta B/B}
\subset
G/B~,$$
where the closure is taken with respect to the (inductive limit) Zariski
topology on $G/B$. In particular,
$$\overline{ I_{\theta,\beta } }\setminus I_{\theta,\beta } =
\{\bar{v}^{-1}\bar{\theta}\frak e, \bar{v}^{-1}\bar{\theta}\bar{r_\beta} \frak
e\}~,$$
where $\bar{v}$ is a preimage of $v$ in $N(T)$.  By Lemma (8.5), it is easy to
see that any closed  $T$-stable subset of $U^{-}\frak e$  (under the induced
subspace topology on $U^{-}\frak e \subset G/B$) contains $\frak e$. Hence  (if
$v\neq \theta$)
$$ I_{\theta,\beta } \subset U^{-} \frak e \Leftrightarrow \frak e \in
\overline{I_{\theta,\beta }} \cap U^{-}\frak e \Leftrightarrow
\bar{v}^{-1}\bar{\theta}\bar{r_\beta} \frak e = \frak e$$
i.e. $v=\theta r_\beta$.   Again by  the Bruhat decomposition for $SL(2)$,  it
is easy to see that in this  case (i.e.  $\theta r_\beta =v$)
$v^{-1}\theta (U_{-\beta}\setminus e)\frak e=  (U_{-\beta}\setminus e)\frak e$.

 But  by \cite{BGG, Corollary 2.3},
 $$
 \{ \beta \in \Delta _+^{\text{re}} : \beta \in \Delta _+
 \cap
 (r_\beta v^{-1}\Delta _-) \text{ and }
 v r_\beta \leq w \}
 =
 \{ \beta \in S(w,v) : v< vr_\beta  \}.
 $$
  This proves the lemma. \qed
 \enddemo

 Now by an argument identical to the proof of Theorem
 (5.5), we obtain the following.

 \proclaim{(8.9) Theorem} Theorem (5.5) is true for an
 arbitrary Kac-Moody group.\qed
 \endproclaim
\flushpar
{\bf (8.10)} {\it Remarks.} Even though we have taken the base field to be the
field $\Bbb C$ of complex numbers, all the results of the paper carry over
(with the same proofs) to an arbitrary algebraically closed field of char. $0$.

Also, by a result of Polo \cite{P, \S4.1}, the dimension of the Zariski tangent
space $Z_\frak v(X_w)$ is independent of the char. of the field. In particular,
a point $\frak v \in X_w$ is smooth in char. $0$ if and only if it is smooth in
any char. $p$. So our smoothness criterion (as in Theorem 5.5(b)) works in
arbitrary char. $p$.

\vskip6ex

 \centerline{{\bf REFERENCES}}
 \roster

\item"[A]" Andersen, H. H., {\it Schubert varieties and Demazure character
formula}, Invent. Math. {\bf 79} (1985), 611-618.

  \item"[BGG]" Bernstein, I. N., Gel'fand, I. M., and
 Gel'fand, S. I.,
 {\it Schubert cells and cohomology of the spaces $G/P$},
 Russian Math. Surveys {\bf 28} (1973), 1--26.

  \item"[Bo]"  Boe, B. D., {\it Kazhdan--Lusztig
 polynomials for hermitian
 symmetric spaces}, Trans. A.M.S. {\bf 309} (1988),
 279--294.

 \item"[B]" Bourbaki, N., {\sl ``Groupes et alg\`{e}bres
 de Lie}, Chap. IV--VI," Hermann, Paris 1968.

\item"[C]" Carrell, J. B., {\it The Bruhat graph of a
 Coxeter group, a
 conjecture of Deodhar, and rational smoothness of
 Schubert varieties}, Proc. Symp. Pure Math. {\bf 56}
 (1994) (edited by W. J. Haboush and B. J. Parshall),
 53--61.

\item"[C2]" Carrell, J. B., {\it On the singular locus of a Schubert
variety:
A survey}, Preprint (1994).

\item"[CPS]" Cline, E., Parshall, B., Scott, L., {\it Cohomology, hyperalgebras
and representations}, J. of Algebra {\bf 63} (1980),
98--123.

 \item"[D]" Deodhar, V. V.,  {\it Local Poincar\'e
 duality and non--singularity of Schubert varieties},
 Comm. in Algebra {\bf 13}, no. 6 (1985), 1379--1388.

\item"[D2]" Deodhar, V. V., {\it A brief survey of Kazhdan-Lusztig theory and
related topics}, Preprint (1993).

  \item"[Dy]" Dyer, M. J.,  {\it The nil Hecke ring and
 Deodhar's
 conjecture on Bruhat intervals}, Invent. Math. {\bf 111}
 (1993), 571--574.

\item"[F]" Fulton, W., {\sl ``Introduction to Intersection Theory in Algebraic
Geometry}," CBMS regional conference series in Mathematics no. {\bf 54} (1984),
Am. Math. Soc..

\item"[Ha]" Harris, J., ``{\sl Algebraic Geometry}," Springer--Verlag,
Berlin--Heidelberg--New York (1992).

   \item"[H]" Hartshorne, R., ``{\sl Algebraic Geometry}, "
 Springer--Verlag, Berlin--Heidelberg--New York (1977).

   \item"[J]" Jantzen, J. C., ``{\sl Moduln mit einem
 hochsten Gewicht},"
 LNM vol. {\bf 750}, Springer--Verlag,
 Berlin--Heidelberg--New York (1979).

 \item"[Jo]" Joseph, A., {\it On the variety of a highest
 weight module},
 J. of Algebra {\bf 88} (1984), 238--278.

 \item"[Jo2]" Joseph, A., {\it On the Demazure character formula}, Ann. Sci.
\,Ec. Norm. Sup. {\bf 18} (1985), 389--419.

 \item"[K]" Kac, V. G., {\sl ``Infinite dimensional Lie
 algebras}," Third edition, Cambridge
 University Press (1990).

\item"[KP]" Kac, V. G.,  and Peterson, D. H.,  {\it Regular functions on
certain infinite dimensional groups}, In: ``Arithmetic and Geometry," (ed. by
M. Artin and J. Tate), Birkhauser (1983), 141-166.

  \item"[KL]" Kazhdan, D., and Lusztig, G., {\it
 Representations of Coxeter groups and Hecke algebras},
 Invent. Math. {\bf 53} (1979), 165--184.

  \item"[KK1]" Kostant, B., and Kumar, S., {\it The nil
 Hecke ring and cohomology of $G/P$ for a Kac--Moody group
 $G$},  Advances in Math. {\bf 62} (1986), 187--237.

  \item"[KK2]" Kostant, B., and Kumar, S., {\it
 $T-$equivariant $K-$theory of generalized flag
 varieties}, J. Diff. Geom. {\bf 32} (1990), 549--603.

 \item"[Ku]" Kumar, S., {\it Demazure character formula
 in arbitrary Kac--Moody setting}, Invent. Math. {\bf 89}
 (1987), 395--423.

 \item"[Ku2]" Kumar, S., {\it The nil Hecke ring and singularity of Schubert
varieties}, In: ``Lie Theory and Geometry (in honor of Bertram Kostant),"  (ed.
by J. -L. Brylinski et. al.), Progress in Math vol. 123, Birkhauser (1994),
497-507.

 \item"[L]" Lakshmibai, V.,  {\it
 Singular loci of Schubert varieties for classical
 groups}, Bull. A.M.S. {\bf 16} (1987), 83--90.

  \item"[LS]" Lakshmibai, V., and Seshadri C. S., {\it
 Singular locus of a Schubert variety
 }, Bull. A.M.S. {\bf 11} (1984), 363--366.

 \item"[Ma]" Mathieu, O., {\it Formules de caract\`eres pour les alg\`ebres de
Kac-Moody g\'en\'erales}, Ast\'erisque {\bf 159-160} (1988), 1-267.

  \item"[M]" Mumford, D., {\sl ``The Red Book of Varieties
 and Schemes},"  LNM vol. {\bf 1358}, Springer--Verlag
 (1988).

  \item"[P]" Polo, P., {\it On Zariski tangent spaces of
 Schubert varieties, and a proof of a conjecture of
 Deodhar}, Preprint (1993).

\item"[Ra]" Ramanathan, A., {\it Schubert varieties are arithmetically
Cohen-Macaulay}, Invent. Math. {\bf 80} (1985), 283-294.

\item"[Ra2]" Ramanathan, A., {\it Equations defining Schubert varieties and
Frobenius splitting of diagonals}, Publ Math. IHES no. {\bf 65} (1987), 61-90.

  \item"[R]" Rossmann, W., {\it Equivariant
 multiplicities on complex
 varieties}, Ast\'{e}risque no. {\bf 173--174} (1989),
 313--330.

  \item"[Ry]"  Ryan, K. M., {\it On Schubert varieties in
 the flag
 manifold of SL($n, \Bbb C$)}, Math. Annalen {\bf 276}
 (1987), 205--224.

\item"[Se]" Seshadri, C. S., {\it Line bundles on Schubert varieties},
In: ``Vector Bundles on Algebraic Varieties," Tata Institute of Fundamental
Research, Bombay (1984), 499--528.

\item"[S]" Slodowy, P., {\it On the geometry of Schubert varieties attached to
Kac-Moody Lie-algebras}, Canad. Math. Soc. Conf. Proc. on `Algebraic Geometry'
(Vancouver) vol. {\bf 6} (1984), 405--442.

 \endroster

 \vskip.5cm
 \flushpar
 Department of Mathematics \flushpar
 University of North Carolina \flushpar
 Chapel Hill,
 N.C. 27599-3250 \flushpar
 U.S.A.

 \enddocument